%%%%%%%%%%%%%%%%%%%%%%%%%%%%%%%%%%%%%%%%%%%%%%%%%%%%%%%%%%
%%%%%%%%%%%%%%%%%%%%%%%%%%%%%%%%%%%%%%%%%%%%%%%%%%%%%%%%%%
%%
%%     This is the AMS-LaTeX file:
%%
%%     Sprekels-Tr\"oltzsch 1
%%       %%     Sparse optimal control of a phase field system with 
%%     singular potentials arising in the  
%%     modeling of tumor growth 
%%     20.08.2020 
%%
%%%%%%%%%%%%%%%%%%%%%%%%%%%%%%%%%%%%%%%%%%%%%%%%%%%%%%%%%%
%%%%%%%%%%%%%%%%%%%%%%%%%%%%%%%%%%%%%%%%%%%%%%%%%%%%%%%%%%

\def\LaTeX{%
  \let\Begin\begin
  \let\End\end
  \let\salta\relax
  \let\finqui\relax
  \let\futuro\relax}

\def\UK{\def\our{our}\let\sz s}
\def\USA{\def\our{or}\let\sz z}

\UK
%\USA

%%%%%%%%%%%%%%%%%%%%%%%%%%%%%%%%%

% scegliere fra \TeX e \LaTeX  e fra  \UK oppure \USA

%\TeX
\LaTeX

%\UK
\USA

%%%%%%%%%%%%%%%%%%%%%%%%%%%%%%%%%
%% page layout
%%%%%%%%%%%%%%%%%%%%%%%%%%%%%%%%%

\salta

\documentclass[twoside,12pt]{article}
\setlength{\textheight}{24cm}
\setlength{\textwidth}{16cm}
\setlength{\oddsidemargin}{2mm}
\setlength{\evensidemargin}{2mm}
\setlength{\topmargin}{-15mm}
\parskip2mm

%%%%%%%%%%%%%%%%%%%%%%%%%%%%%%%%%
%% packages
%%%%%%%%%%%%%%%%%%%%%%%%%%%%%%%%%

%\usepackage{color}
\usepackage[usenames,dvipsnames]{color}
\usepackage{amsmath}
\usepackage{amsthm}
\usepackage{amssymb}
\usepackage[mathcal]{euscript}
\usepackage{cite}

\usepackage[ulem=normalem,draft]{changes}
%
%		COLORS FOR CORRECTIONS
%
% do the same, please (i.e., don't use the standard {\color{red} text} or similar): 
% just choose the color you prefer in \def\yourname

% EXAMPLE OF USE:  \fredi{I want this to become blue}
%
%IF YOU LATER WANT TO LET THE COLOR DISAPPEAR, ACTIVATE \def\fredi #1{{#1}} BELOW
 
\definecolor{viola}{rgb}{0.3,0,0.7}
\definecolor{ciclamino}{rgb}{0.5,0,0.5}
\definecolor{rosso}{rgb}{0.85,0,0}

\def\juerg #1{{\color{red}#1}}

\def\neuj #1{{\color{green}#1}}
\def\juerg #1{#1}

\def\neuj #1{{#1}}

%%%%%%%%%%%%%%%%%%%%%%%%%%%%%%%%%
%% you may adjust the baseline
%%%%%%%%%%%%%%%%%%%%%%%%%%%%%%%%%

%\renewcommand{\baselinestretch}{0.975}

%%%%%%%%%%%%%%%%%%%%%%%%%%%%%%%%%
%% bibliographystyle
%%%%%%%%%%%%%%%%%%%%%%%%%%%%%%%%%

\bibliographystyle{plain}

%%%%%%%%%%%%%%%%%%%%%%%%%%%%%%%%%
%% environments
%%%%%%%%%%%%%%%%%%%%%%%%%%%%%%%%%

%
\newtheorem{theorem}{Theorem}[section]

\newtheorem{lemma}[theorem]{Lemma}

\finqui

\def\Bthm{\Begin{theorem}}
\def\Ethm{\End{theorem}}
\def\Blem{\Begin{lemma}}
\def\Elem{\End{lemma}}

\def\Brem{\Begin{remark}\rm}
\def\Erem{\End{remark}}

\def\Bdim{\Begin{proof}}
\def\Edim{\End{proof}}
\def\Bcenter{\Begin{center}}
\def\Ecenter{\End{center}}
\let\non\nonumber

%%%%%%%%%%%%%%%%%%%%%%%%%%%%%%%%%
%% macros
%%%%%%%%%%%%%%%%%%%%%%%%%%%%%%%%%

% macro salvate

% sottosezioni non numerate

\def\step #1 \par{\medskip\noindent{\bf #1.}\quad}
\def\jstep #1: \par {\vspace{2mm}\noindent\underline{\sc #1 :}\par\nobreak\vspace{1mm}\noindent}

% versioni inglesi (UK) o americane (USA)

% bold, cal e mathop

\def\multibold #1{\def\arg{#1}%
  \ifx\arg\pto \let\next\relax
  \else
  \def\next{\expandafter
    \def\csname #1#1#1\endcsname{{\bf #1}}%
    \multibold}%
  \fi \next}

\def\pto{.}

\def\multical #1{\def\arg{#1}%
  \ifx\arg\pto \let\next\relax
  \else
  \def\next{\expandafter
    \def\csname cal#1\endcsname{{\cal #1}}%
    \multical}%
  \fi \next}

% operatori

\def\multimathop #1 {\def\arg{#1}%
  \ifx\arg\pto \let\next\relax
  \else
  \def\next{\expandafter
    \def\csname #1\endcsname{\mathop{\rm #1}\nolimits}%
    \multimathop}%
  \fi \next}

\multibold
qwertyuiopasdfghjklzxcvbnmQWERTYUIOPASDFGHJKLZXCVBNM.

\multical
QWERTYUIOPASDFGHJKLZXCVBNM.

\multimathop
diag dist div dom mean meas sign supp .

% accorpamenti di formule citate:
% uso  \accorpa {prima}{seconda}
%      \Accorpa\cs prima seconda (con il comodo blank anche dopo)
% NB: \Accorpa definisce \cs come l'accorpamento delle due citazioni
% e scrive sul file.log

\def\Accorpa #1#2 #3 {\gdef #1{\eqref{#2}--\eqref{#3}}%
  \wlog{}\wlog{\string #1 -> #2 - #3}\wlog{}}

% macro comode

\def\<#1>{\mathopen\langle #1\mathclose\rangle}

\def\intQ{\int_Q}
\def\iO{\int_\Omega}

\def\dt{\partial_t}
\def\dn{\partial_{\bf n}}
\def\S{{\cal S}}

\def\X{{\cal X}}
\def\Y{{\cal Y}}
\def\Uh{{\cal U}}

\def\checkmmode #1{\relax\ifmmode\hbox{#1}\else{#1}\fi}

% insiemi numerici

\def\bu{{\bf u}}

\def\bl{{\boldsymbol \lambda}}

\def\erre{{\mathbb{R}}}
\def\rz{{\mathbb{R}}}

\def\enne{{\mathbb{N}}}

% spazi di funzioni a valori vettoriali su [0,T], [0,t], [0,s], [0,+\infty), [\delta,T]

% Come ricordare: in generale i simboli L H W  C da soli per gli spazi su (0,T)
% gli stessi raddoppiati per (0,+\infty)
% aggiunta di t o s al simbolo per (0,t) e (0,s)
% aggiunta di d al simbolo semplice o doppio per intervalli (\delta,T) e (\delta,+\infty)
% il simbolo C e i suoi derivati mettono le quadre anziche' le tonde

% Esempi   \L2V   \L\infty\Vp   \W{1,1}H   \C0H   \LL2V   \CC0\Vp   \Ld2V  \CCdH

%\def\LL {\spazioinf L}
%\def\HH {\spazioinf H}
%\def\WW {\spazioinf W}
%\def\CC #1#2{C^{#1}([0,+\infty);#2)}

% spazi di funzioni su \Omega, \Gamma, Q e \Sigma

\def\Lx #1{L^{#1}(\Omega)}
\def\Hx #1{H^{#1}(\Omega)}

\def\Ldue{\Lx 2}

\def\Huno{\Hx 1}
\def\Hdue{\Hx 2}

\def\Liq{{L^\infty(Q)}}

% spazi di funzioni su Q e S

% lettere greche

%\let\badtheta\theta
%\let\theta\vartheta

\let\phi\varphi
\let\vp\varphi

\def\opi{{\overline\psi}}

\let\TeXchi\chi                         % new \chi, exactly on the baseline
\newbox\chibox
\setbox0 \hbox{\mathsurround0pt $\TeXchi$}
\setbox\chibox \hbox{\raise\dp0 \box 0 }
\def\chi{\copy\chibox}

% quadratino di fine dimostrazione

% abbreviazioni specifiche del lavoro

\def\ubar{\overline{\bf u}}
\def\uebar{\overline{u}_1}
\def\uzbar{\overline{u}_2}

\let\hat\widehat

\def\uad{{\cal U}_{\rm ad}}
\def\bmu{\overline\mu}
\def\bvp{\overline\varphi}
\def\bsigma{\overline\sigma}

%\def\HH{\hat\calH}
%\def\VV{\hat\calV}

%%%%%%%%%%%%%%%%%%%%%%%%%%%%%%
\Begin{document}
%%%%%%%%%%%%%%%%%%%%%%%%%%%%%%%%%

%%%%%%%%%%%%%%%%%%%%%%%%%%%%%%%%%
%% front page
%%%%%%%%%%%%%%%%%%%%%%%%%%%%%%%%%

%
\title{Sparse optimal control of a phase field\\
system with singular potentials arising \\in the
     modeling of tumor growth}
\author{}
\date{}
\maketitle
\Bcenter
\vskip-1.5cm
{\large\sc J\"urgen Sprekels$^{(1)}$}\\
{\normalsize e-mail: {\tt sprekels@wias-berlin.de}}\\[.25cm]
{\large\sc Fredi Tr\"oltzsch$^{(2)}$}\\
{\normalsize e-mail: {\tt troeltz@math.tu-berlin.de}}\\[.5cm]
$^{(1)}$
{\small Department of Mathematics}\\
{\small Humboldt-Universit\"at zu Berlin}\\
{\small Unter den Linden 6, D-10099 Berlin, Germany}\\[2mm]
{\small and}\\[2mm]
{\small Weierstrass Institute for Applied Analysis and Stochastics}\\
{\small Mohrenstrasse 39, D-10117 Berlin, Germany}\\[5mm]
$^{(2)}$
{\small Institut f\"ur Mathematik}\\
{\small Technische Universit\"at Berlin}\\
{\small Strasse des 17. Juni 136, D-10623 Berlin, Germany}
\Ecenter
\Begin{abstract}
\noindent In this paper, we study an optimal control problem for a nonlinear system of reaction-diffusion
equations that constitutes a simplified and relaxed version of a thermodynamically
 consistent phase field model for tumor growth originally introduced in \cite{GLSS}. The model takes the effect of chemotaxis into account but neglects velocity contributions. {The unknown quantities of the governing state equations are the chemical potential, the (normalized) tumor fraction, and the nutrient extra-cellular water concentration. The equation governing the evolution of the tumor fraction is dominated by the variational derivative of a double-well potential which may be of singular (e.g., logarithmic) type.} In contrast to the recent paper \cite{CSS1} on the same system, we consider in this paper sparsity effects, which means that the cost functional contains a nondifferentiable (but convex) contribution like the $L^1-$norm. For such problems, we derive first-order necessary optimality conditions and conditions for directional sparsity, both with respect to space and time, where the latter case is of particular interest for practical medical applications in
  which the control variables are given by the administration of cytotoxic drugs or by the supply of nutrients. In addition to these results, we prove that the corresponding control-to-state operator is twice continuously differentiable between suitable Banach spaces, using the implicit function theorem. This result, which complements and sharpens a differentiability result derived in \cite{CSS1}, constitutes a prerequisite for a future derivation of second-order sufficient optimality conditions.


\End{abstract}
\salta
\pagestyle{myheadings}
\newcommand\testopari{\sc Sprekels --- Tr\"oltzsch}
\newcommand\testodispari{\sc Sparse optimal control for phase field tumor growth model}
\markboth{\testopari}{\testodispari}
\finqui
%
%%%%%%%%%%%%%%%%%%%%%%%%%%%%%%%%%
%% very beginning
%%%%%%%%%%%%%%%%%%%%%%%%%%%%%%%%%

\section{Introduction}
\label{INTRO}
\setcounter{equation}{0}

Let $\Omega\subset\erre^3$ denote some open, bounded and connected set having a smooth boundary $\Gamma=\partial\Omega$
and unit outward normal $\,{\bf n}$.
We denote by $\dn$ the outward normal derivative to $\Gamma$.  Moreover, we fix some final time $T>0$ and
introduce for every $t\in (0,T]$ the sets $Q_t:=\Omega\times (0,t)$ and $\Sigma_t:=\Gamma\times (0,t)$,
where we put, for the sake of brevity, $Q:=Q_T$ and $\Sigma:=\Sigma_T$. We then consider the following optimal control
problem:

\vspace{3mm}\noindent
(${\cal CP}$) \quad Minimize the cost functional
\begin{align} 
\label{cost} 
{\cal J}((\mu,\vp,\sigma),{\bf u})
   :=&\,  \frac{\beta_1}2 \intQ |\vp-\widehat \vp_Q|^2
  + \frac{\beta_2}2 \iO |\vp(T)-\widehat\vp_\Omega|^2
  + \frac{\nu}2 \intQ |{\bf u}|^2 \,+\,\kappa\,g({\bf u})\,,
\end{align} 
subject to  the state system  
\begin{align}
\label{ss1}
&\alpha\dt\mu+\dt\phi-\Delta\mu=(P\sigma-A-u_1)h(\vp) \quad\mbox{in }\,Q\,,\\[1mm]
\label{ss2}
&\beta\dt\vp-\Delta\vp+F'(\vp)=\mu+\chi\,\sigma \quad\mbox{in }\,Q\,,\\[1mm]
\label{ss3}
&\dt\sigma-\Delta\sigma=-\chi\Delta\vp+B(\sigma_s-\sigma)-E\,\sigma h(\vp)+u_2 \quad\mbox{in }\,Q\,,\\[1mm]
\label{ss4}
&\dn \mu=\dn\vp=\dn\sigma=0 \quad\mbox{on }\,\Sigma\,,\\[1mm]
\label{ss5}
&\mu(0)=\mu_0,\quad \vp(0)=\vp_0,\quad \sigma(0)=\sigma_0\quad \mbox{in }\,\Omega\,,
\end{align}
and to the control constraint
\begin{equation}
\label{uad}
{\bf u}=(u_1,u_2)\in\uad\,.
\end{equation}
Here, the constants $\beta_1,\beta_2$ are nonnegative, while $\,\nu\,$ and $\,\kappa\,$ are positive. Moreover,
 $\widehat\vp_Q$ and $\widehat \vp_\Omega$ are given target functions, and $g:{\cal U}\to [0,+\infty)$ is a 
 nonnegative and 
convex, but not necessarily differentiable, functional on the control space  
\begin{equation}
\label{defU}
\Uh:= L^\infty(Q)^2.
\end{equation}
Moreover, $\uad$ is a suitable bounded, closed and convex subset of ${\cal U}$. Since we are interested in sparse 
controls in this note, typical (nondifferentiable) examples for the functional $g$ are given by
\begin{align}
\label{sparse}
&g(\bu)=\|\bu\|_{L^1(Q)}=\intQ |\bu(x,t)|\,dx\,dt,\\
\label{dtsparse}
&g(\bu)={\int_0^T\Bigl(\iO |\bu(x,t)|^2\,dx\Bigr)^{1/2} dt,}\\
&g(\bu)={\iO\Bigl(\int_0^T|\bu(x,t)|^2\,dt\Bigr)^{1/2} dx.}
\label{dxsparse}
\end{align}
The functionals in \eqref{dtsparse} and \eqref{dxsparse} are associated with the notion of {\em {directional}
sparsity} (with respect to $\,t\,$ and to $\,x$, respectively). Since we have two control variables in our 
system, we could ``mix'' the sparsity directions by taking different ones for $u_1$ and $u_2$; also, different weights
could be given to the directions. For the sake of avoiding unnecessary technicalities, we restrict ourselves to 
the simplest case here.  

 The state system \eqref{ss1}--\eqref{ss5} constitutes a simplified and relaxed version of a thermodynamically
 consistent phase field model for tumor growth that includes the effect of chemotaxis and was 
 originally introduced in \cite{GLSS}. Indeed, the velocity
 contributions in \cite{GLSS} were neglected, and the two relaxation terms $\,\alpha\dt\mu\,$ are $\,\beta\dt\vp\,$ have been
 added. We note that a different thermodynamically consistent model was introduced in \cite{HZO} and 
 studied mathematically in \cite{CGH, CGRS1,CGRS2,CGRS3}, where \cite{CGRS3} focused on optimal control
 problems. In this connection, we also refer to \cite{CGS24}.
 
 In all of the abovementioned models, the unknowns $\,\mu,\vp,\sigma\,$ stand for the chemical potential, the normalized
 tumor fraction, and the nutrient extra-cellular water concentration, in this order. The quantity $\sigma$ is
usually normalized between $\,0\,$ and $\,1$, where these values model nutrient-poor and nutrient-rich cases.
The variable $\,\vp\,$ plays the role of an order parameter and is usually taken between the values $\,-1\,$ and $\,+1$, 
which represent the healthy cell case and the tumor phase, respectively. The capital letters $\,A,B,{E},P,\chi\,$
denote positive coefficients that stand for the apoptosis rate, nutrient supply rate, nutrient consumption
 rate, and chemotaxis coefficient, in this order. In addition, let us point out that the contributions
$\,\chi\,\sigma\,$ and $\,\chi\,\Delta\vp\,$ model pure chemotaxis. Furthermore, the nonlinear function
$\,h\,$ has been considered in \cite{GLSS} as an interpolation function satisfying $h(-1)=0$ and $h(1)=1$, so 
that the mechanisms modeled by the terms $\,(P\sigma-A-u_1)h(\vp)\,$ and $\,{E}\sigma h(\vp)\,$ are switched off in the 
healthy tissue (which corresponds to $\vp=-1$) and are fully active in the tumorous case $\vp=1$. Moreover, the 
term $\,\sigma_s\,$ is a nonnegative constant that models the nutrient concentration in a pre-existing 
vasculature. 
 
Very important for the evolution of the state system is the nonlinearity $F$, which is  assumed 
to be a double-well potential. Typical examples are given by the regular and logarithmic potentials, which are
given, in this order, by
\begin{align}
\label{regpot}
&F_{\rm reg}(r)=\frac 14 \left(1-r^2\right)^2 \quad\mbox{for $\,r\in\rz$,}\\[1mm]
\label{logpot}
&F_{\rm log}(r)=(1+r)\,\ln(1+r)+(1-r)\,\ln(1-r)-k r^2\quad\mbox{for $\,r\in (-1,1)$,}
\end{align}
where $k>1$ so that $F_{\rm log}$ is nonconvex. Observe that $F_{\rm log}$ is very relevant  in the 
applications, where $F'_{\rm log}(r)$ becomes unbounded as $r\to\pm 1$.

In this paper, we work with two source controls that act in the phase equation and in the nutrient equation,
respectively. The control variable $u_1$ in the phase equation models the application of a cytotoxic drug into 
the system; it is multiplied by $h(\vp)$ in order to have the action
only in the spatial region where the tumor cells are located. On the other hand, the control $u_2$ can model  
either an external medication or some nutrient supply. In this connection, sparsity of the control is
highly desirable: indeed, if a distributed cytotoxic drug is to be administered, this should be done only
where it does not harm healthy tissue, which calls for {directional} sparsity with respect to space;
on the other hand, and even more importantly, cytotoxic drugs should only be applied for very short periods of time,
in order to prevent the tumor cells from developing a resistance against the drug. This, of course, calls for a 
{directional}
sparsity with respect to time.

Optimal control problems for the system \eqref{ss1}--\eqref{ss5} have recently been studied in \cite{CSS1}, where the
cost functional, while containing some additional quadratic terms, did not have a nondifferentiable 
contribution, i.\,e., we had $g\equiv 0$. However, besides existence of optimal controls, it was shown in 
\cite{CSS1} that the control-to-state operator is Fr\'echet differentiable between suitable function spaces, and
first-order necessary optimality conditions in terms of the adjoint state variables were derived. The 
Fr\'echet differentiability was shown ``directly'' without using the implicit function theorem, and therefore
the existence of higher-order derivatives was  not proved. Note that the existence of second-order derivatives
forms a prerequisite for deriving second-order sufficient optimality conditions and efficient numerical 
techniques. To pave the way for such an analysis (which shall not be given in this paper), we 
have decided to include a proof of the Fr\'echet differentiability of the control-to-state 
operator via the implicit function theorem.  

Another novelty of this paper is the discussion of optimal controls that are sparse with respect to the time. Since the seminal paper \cite{stadler2009}, sparse optimal controls have been discussed extensively in the literature. Directional sparsity was introduced in \cite{herzog_stadler_wachsmuth2012, herzog_obermeier_wachsmuth2015} and extended to semilinear parabolic optimal control problems in
\cite{casas_herzog_wachsmuth2017}. Sparse optimal controls for reaction-diffusion equations  were investigated in \cite{casas_ryll_troeltzsch2013,casas_ryll_troeltzsch2015}. 

Although the main techniques of the analysis for sparse controls are known from the abovementioned papers, 
a discussion of sparsity for the control of the system of reaction-diffusion equations \eqref{ss1}--\eqref{ss5} seems to be worth investigating in view of its medical background. In this connection,  temporal sparsity is particularly interesting. It means that drugs are not needed in certain time periods. For the control of the class of reaction-diffusion equations \eqref{ss1}--\eqref{ss5}, the investigation of sparse controls is new.

The paper is organized as follows: in the subsequent Section 2, we give the general setting of the problem, 
and recall known  well-posedness results for the state system \eqref{ss1}--\eqref{ss5}. Moreover, we employ
the implicit function theorem to show that the control-to-state is twice continuously Fr\'echet differentiable between suitable Banach spaces, thereby sharpening the differentiability result of \cite{CSS1}. Section 3 then deals with
first-order necessary optimality conditions for the problem $({\cal CP})$, and the final Section 4 brings a 
discussion of the sparsity of optimal controls.
  
Throughout this paper, we make repeated use of H\"older's inequality, of the elementary Young's inequality
\begin{equation}
\label{Young}
a\,b\,\le\,\gamma |a|^2\,+\,\frac 1{4\gamma}\,|b|^2\quad\forall\,a,b\in\erre, \quad\forall\,\gamma>0,
\end{equation}
as well as the continuity of the embeddings $H^1(\Omega)\subset L^p(\Omega)$ for $1\le p\le 6$ and 
$\Hdue\subset C^0(\overline\Omega)$. Notice that the latter embedding is also compact, while this holds true
for the former embeddings only if $p<6$. 
%We will also use the denotations
%\begin{align}
%\label{defQt}
%&Q^t:=\Omega\times (t,T),\quad \Sigma^t:=\Gamma\times (t,T),\quad\mbox{for }\,0\le t<T.
%\end{align}
Moreover, throughout the paper, for a Banach space $\,X\,$ we denote by $\,\|\,\cdot\,\|_X\,$
the norm in the space $X$ or in a power of it, and by $\,X^*\,$ its dual space. 
 The only exemption from this rule applies to the norms of the
$\,L^p\,$ spaces and of their powers, which we often denote by $\|\,\cdot\,\|_p$, for
$\,1\le p\le +\infty$. As usual, for Banach spaces $\,X\,$ and $\,Y\,$ we introduce the linear space
$\,X\cap Y\,$ which becomes a Banach space when equipped with its natural norm $\,\|u\|_{X \cap Y}:=
\|u\|_X\,+\,\|u\|_Y$, for $\,u\in X\cap Y$.

%%%%%%%%%%%%%%%%%%%%%%%%%%%%%%%%%%%%%%%%%%%%%%%%%%%%%%%%%%%%%%%%%%%%%%

\section{General setting and properties of the state system}
\label{STATE}
\setcounter{equation}{0}

In this section, we introduce the general setting of our control 
problem and state some results on the state system \eqref{ss1}--\eqref{ss5}. To begin with, we
recall the definition~\eqref{defU} of~${\cal U}$ and introduce the spaces
\begin{align}
  & H := \Ldue \,, \quad  
  V := \Huno\,,   \quad
  W_{0} := \{v\in\Hdue:\dn v=0 \,\mbox{ on $\,\Gamma$}\}.
  \label{defHVW}
\end{align}
By $\,(\,\cdot\,,\,\cdot\,)$ we denote the standard inner product in $\,H$.

For the potential $\,F$, we generally assume:

\vspace{2mm}\noindent
{\bf (F1)}\quad $F=F_1+F_2$, where $\,F_1:\erre\to [0,+\infty]\,$ is convex and lower semicontinuous with
\hspace*{11mm}
$\,F_1(0)=0$.

\vspace{1mm}\noindent
{\bf (F2)}\quad There exists an interval $\,(r_-,r_+)\,$ with $\,-\infty\le r_-<0<r_+\le +\infty\,$ such that
the \hspace*{12.2mm} restriction of $F_1$ to $\,(r_-,r_+)\,$ belongs to $\,C^5(r_-,r_+)$.

\vspace{1mm}\noindent
{\bf (F3)}\quad $F_2\in C^5(\erre)$, and $\,F_2'\,$ is globally Lipschitz continuous on $\erre$.

\vspace{1mm}\noindent
{\bf (F4)}\quad It holds $\,\lim_{r\to r_{\pm}} F'(r)=\pm\infty$.                                      

\noindent It is worth noting that both \eqref{regpot} and \eqref{logpot} fit into this framework with the choices
$\,(r_-,r_+)=\erre\,$ and $\,(r_-,r_+)=(-1,1)$, respectively, where in the latter case we extend $F_{\rm log}$
by $\,F_{\rm log}(\pm 1)=2\ln(2)-k\,$ and $\,F_{\rm log}(r)=+\infty\,$ for $\,r\not\in [-1,1]$.

For the initial data, we make the following assumptions:

\vspace{2mm}\noindent
{\bf (A1)}\quad $\vp_0,\mu_0,\sigma_0\in W_{0}$,
% \cap W^{\frac 53,6}(\Omega)$, 
\,\,and\,\, $r_-<\min_{x\in\overline{\Omega}}\,\vp_0(x) \le
\max_{x\in\overline{\Omega}}\,\vp_0(x)<r_+$.

\noindent Notice that {\bf (A1)} entails that $F^{(i)}(\vp_0)\in C^0(\overline\Omega)$,
for $i=0,...,5$.
This condition can be restrictive for the case of singular potentials; for instance, in the case of
the logarithmic potential $\,F_{\rm log}\,$ we have $r_\pm=\pm 1$, so that {\bf (A1)} excludes the
pure phases (tumor and healthy tissue) as initial data.  

For the other data and the target functions, we postulate:

\vspace{2mm}\noindent
{\bf (A2)}\quad $h\in C^3(\erre)\cap W^{3,\infty}(\erre)$, and $\,h\,$ is positive on $(r_-,r_+)$. 

\vspace{1mm}\noindent
{\bf (A3)}\quad $\alpha,\beta,\chi$ are positive constants, while $P,A,B,E,\sigma_s$ are nonnegative constants.

\vspace{1mm}\noindent
{\bf (A4)}\quad $\beta_1,\beta_2$ are nonnegative, and $\nu,\kappa$ are positive.

\vspace{1mm}\noindent
{\bf (A5)}\quad $\vp_Q\in L^2(Q)$ and $\,\vp_\Omega\in \Ldue$.

\vspace{2mm}\noindent
Observe that {\bf (A2)} entails that $h,h',h''$ are Lipschitz continuous on $\erre$. We now assume for the
set of admissible controls:

\vspace{1mm}\noindent
{\bf (A6)}\quad $\uad=\left\{\bu=(u_1,u_2)\in {\cal U}: \underline u _i\le u_i\le\hat u_i\mbox{\, a.e. in 
$Q$}, \,\,\,i=1,2\right\},$\\[1mm]
\hspace*{14mm}where \,$\underline u_i, \hat u_i\in L^\infty(Q)$\, satisfy \,$\underline u_i\le
\hat u_i$\,  a.e. in $\,Q$, $i=1,2$. 

\vspace{2mm}\noindent
Notice that $\uad$ is a nonempty, closed and convex subset of $\,{\cal U}=L^\infty(Q)^2$. In the 
following, it will sometimes be convenient to work with a bounded open superset of $\uad$.
We therefore once and for all fix some $R>0$ such that
\begin{equation}\label{defUR}
{\cal U}_R:=\left\{\bu =(u_1,u_2) \in L^\infty(Q)^2:\,\|\bu\|_\infty<R\right\} \supset \uad.
\end{equation}

The following result concerning the wellposedness of the state system has been shown in \cite[Thm.~2.2]{CSS1}.
\Bthm
Suppose that the conditions {\bf (F1)}--{\bf (F4)}, {\bf (A1)}--{\bf (A3)}, {\bf (A6)}, and \eqref{defUR} are fulfilled.
Then the state system \eqref{ss1}--\eqref{ss5} has for every $\bu=(u_1,u_2)\in {\cal U}_R$ a unique solution
$(\mu,\vp,\sigma)$ with the regularity
\begin{align}
\label{regmu}
&\mu\in  H^1(0,T;H) \cap C^0([0,T];V) \cap L^2(0,T;W_0)\cap  C^0(\overline Q),\\[1mm]
\label{regphi}
&\phi\in W^{1,\infty}(0,T;H)\cap H^1(0,T;V)\cap L^\infty(0,T;W_0) \cap C^0(\overline Q),
\\[1mm]
\label{regsigma}
&\sigma\in H^1(0,T;H)\cap C^0([0,T];V)\cap L^2(0,T;W_0)\cap C^0(\overline Q).
\end{align}
Moreover, there exists a constant $K_1>0$, which depends on $\Omega,T,R,\alpha,\beta$ and the data of the 
system, but not on the choice of $\bu\in{\cal U}_R$, such that
\begin{align}\label{ssbound1}
&\|\phi\|_{W^{1,\infty}(0,T;H)\cap H^1(0,T;V) \cap L^{\infty}(0,T;W_0)\cap C^0(\overline Q)}
\nonumber\\[1mm]
&+\,\|\mu\|_{H^1(0,T;H) \cap C^0([0,T];V) \cap L^2(0,T;W_0)\cap C^0(\overline  (Q))}\nonumber\\[1mm]
&+\,\|\sigma\|_{H^1(0,T;H) \cap C^0([0,T];V) \cap L^2(0,T;W_0)\cap C^0(\overline Q)}\,\le\,K_1\,.
\end{align}
In addition, there are constants $r_*,r^*$, which depend on $\Omega,T,R,\alpha,\beta$ and the data of the 
system, but not on the choice of $\bu\in{\cal U}_R$, such that
\begin{equation}\label{ssbound2}
r_-  <r_*\le\vp(x,t)\le r^*<r_+ \quad\mbox{for all $(x,t)\in \overline Q$}.
\end{equation}
Finally, there is some constant $K_2>0$, which depends on $\Omega,T,R,\alpha,\beta$ and the data of the 
system, but not on the choice of $\bu\in{\cal U}_R$, such that
\begin{equation}
\label{ssbound3}
\max_{i=0,1,2,3}\,\left\|h^{(i)}(\vp)\right\|_{C^0(\overline Q)}\,+\,\max_{i=0,1,2,3,4,5}\,\left\|
F^{(i)}(\vp)\right\|_{C^0(\overline Q)} \,\le\,K_2\,.
\end{equation}
\Ethm
\Brem
{In the original statement of \cite[Thm.~2.2]{CSS1} it was only asserted that $\mu,\varphi,\sigma\in L^\infty(Q)$. Note,
however, that under our assumptions the initial values are smoother than in \cite{CSS1}. In particular, they belong to
$C^0(\overline\Omega)$. Moreover, since
$\varphi\in W^{1,\infty}(0,T;H)\cap L^\infty(0,T;W_0)$, and as $W_0\subset\Hdue$ is compactly embedded in 
$C^0(\overline\Omega)$, it follows from \cite[Sect.~8, Cor.~4]{Simon} that 
$\varphi\in C^0(\overline Q)$ with the corresponding estimate.}

\juerg{
In addition, we have $\dt\varphi\in (L^\infty(0,T;H)\cap L^2(0,T;V)) \subset L^{10/3}(Q)$. Hence, by bringing 
$\,\dt\phi\,$ to the right-hand side of \eqref{ss1}, we see that $\,\mu\,$ solves a standard linear
parabolic Neumann problem whose right-hand side is bounded in $L^{10/3}(Q)$. Now observe that $\frac{10}3>\frac 52$. It 
thus follows from known parabolic regularity results (see, e.g., \cite[Lem.~7.12]{Fredibuch}) that also 
$\mu \in C^0(\overline Q)$ with the corresponding estimate. Finally, we write \eqref{ss3} in the form}
\juerg{
\begin{equation*}
\dt(\sigma-\chi\vp) - \Delta(\sigma-\chi\vp)=-\chi\dt\vp+B(\sigma_s-\sigma)-E\sigma h(\vp)+u_2,
\end{equation*}}
\juerg{
so that $\,z:=\sigma-\chi\vp$ solves a standard linear parabolic Neumann problem whose right-hand side is
bounded in $L^{10/3}(Q)$. Again, \cite[Lem.~7.12]{Fredibuch} 
yields that $z\in C^0(\overline Q)$, and thus $\sigma\in C^0(\overline Q)$, with corresponding estimates 
for the norms.} 
\Erem

\Brem
The {\em separation condition} \eqref{ssbound2} is particularly important for the case of singular 
potentials such as $F_{\rm log}$. Indeed, it guarantees that the phase variable always stays away
from the critical values $\,r_-,r_+$ that usually correspond to the pure phases. In this way, the singularity 
is no longer an obstacle for the analysis; however, the case of pure phases is then excluded, which is
not desirable from the viewpoint of medical applications.
\Erem                      
 
%\vspace{2mm}
Owing to Theorem 2.1, the control-to-state operator $\,\S:\bu=(u_1,u_2)\mapsto 
(\mu,\phi,\sigma)\,$ is well defined as a
mapping between ${\cal U}=L^\infty(Q)^2$ and the Banach space specified by the regularity results \eqref{regmu}--\eqref{regsigma}.
 
We now discuss the Fr\'echet differentiability of  $\S$, considered 
as a mapping between suitable Banach spaces. We remark that in \cite[Thm.~2.6]{CSS1} 
Fr\'echet differentiability was established between $L^2(Q)^2$ and 
$\left(C^0([0,T];H)\cap L^2(0,T;V)\right) \times (H^1(0,T;H)\cap L^\infty(0,T;V))
\times \left(C^0([0,T];H)\cap L^2(0,T;V)\right)$. The proof was a direct one that did not use 
the implicit function theorem. The result was strong enough to derive meaningful first-order necessary
conditions, but it did not admit the derivation of second-order sufficient conditions, since these
require the control-to-state operator to be twice continuously Fr\'echet differentiable. To show such a result,
it is more favorable to employ the implicit function theorem, because, if applicable, it yields
that the control-to-state operator automatically inherits the differentiability order from that of
the involved nonlinearities.
For this, some functional analytic preparations are in order. We first define the linear spaces
\begin{align}
\X\,&:=\,X\times \juerg{ \widetilde X}\times X, \mbox{\,\,\, where }\nonumber\\
X\,&:=\,H^1(0,T;H)\cap C([0,T];V)\cap L^2(0,T;W_0)\cap \neuj{C^0(\overline Q)},\nonumber\\
\juerg{\widetilde X}\,&\neuj{:=\,W^{1,\infty}(0,T;H)\cap H^1(0,T;V)\cap L^\infty(0,T;W_0)\cap C^0(\overline Q)}, 
\label{calX}
\end{align}
which are Banach spaces when endowed with their natural norms. Next, we introduce the linear space
\begin{align}
&\Y\,:=\,\bigl\{(\mu,\phi,\sigma)\in \calX: \,\alpha\dt\mu+\dt\phi-\Delta\mu\in \Liq, 
\,\,\,\beta\dt\phi-\Delta\phi-\mu\in \Liq,\nonumber\\
& \hspace*{14mm}\dt\sigma-\Delta\sigma+\chi\Delta\phi\in\Liq\bigr\},
\label{calY}
\end{align}
which becomes a Banach space {when endowed} with the norm
\begin{align}
\|(\mu,\phi,\sigma)\|_{\Y}\,:=\,&\|(\mu,\phi,\sigma)\|_\calX
\,+\,\|\alpha\dt\mu+\dt\phi-\Delta\mu\|_{\Liq}
\,+\,\|\beta\dt\phi-\Delta\phi-\mu\|_{\Liq}\nonumber\\
&+\,\|\dt\sigma-\Delta\sigma+\chi\Delta\phi\|_{\Liq}\,.
\label{normY}
\end{align}

Finally, we fix constants $\tilde r_-,\tilde r_+$ such that
\begin{equation}
\label{openphi}
r_-<\tilde r_-<r_*\le r^*<\tilde r_+<r_+,
\end{equation}  
with the constants introduced in {\bf (F2)} and \eqref{ssbound2}. We then consider the set
\begin{equation}
\label{Phi}
\Phi\,:=\,\bigl\{\phi,\mu, \sigma)\in \Y: \tilde r_-<\phi(x,t) <\tilde r_+ \,
\mbox{ for all }\,(x,t)\in\overline Q\bigr\},
\end{equation}
which is obviously an open subset of the space $\Y$.  

We first prove an auxiliary result for the linear initial-boundary value problem
\begin{align}\label{aux1}
&\alpha\dt\mu+\dt\phi-\Delta\mu\,=\,\lambda_1\left[P\sigma h(\bvp)+(P\bsigma-A-\uebar)h'(\bvp)\phi\right]-\lambda_2
k_1 h(\bvp)+\lambda_3 f_1\quad\mbox{in }\,Q,\\[0.5mm]
\label{aux2}
&\beta\dt\phi-\Delta\phi-\mu\,=\,\lambda_1\left[\chi\,\sigma-F''(\bvp)\phi\right]+\lambda_3 f_2 \quad\mbox{ in \,$Q$},
\\[0.5mm]
\label{aux3}
&\dt\sigma-\Delta\sigma+\chi\Delta\phi\,=\,\lambda_1\left[-B\sigma-{E}\sigma h(\bvp)-{E}\bsigma h'(\bvp)\phi\right]
+\lambda_2 k_2+\lambda_3 f_3\quad\mbox{ in \,$Q$},\\[0.5mm]
\label{aux4}
&\dn\mu\,=\,\dn\phi\,=\,\dn\sigma\,=\,0 \quad\mbox{ on \,$\Sigma$},\\[0.5mm]
\label{aux5}
&\mu(0)\,=\,\lambda_4 \mu_0,\quad \phi(0)\,=\,\lambda_4\phi_0, \quad \sigma(0)\,=\,\lambda_4
\sigma_0,\,\,\,\mbox{ in }\,\Omega,
\end{align}
which for $\lambda_1 = \lambda_2 = 1$ and $\lambda_3 = \lambda_4 = 0$ coincides  with the linearization of the state equation at $((\uebar,\uzbar),(\bmu,\bvp,\bsigma))$. We will need this slightly more general version later for the application of the implicit function theorem. 
\Blem \label{L2.4}
Suppose that $\lambda_1,\lambda_2,\lambda_3,\lambda_4 \in\{0,1\}$ are given
and that the assumptions {\bf (F1)}--{\bf (F4)},
{\bf (A1)}--{\bf (A3)}, {\bf (A6)}, and \eqref{defUR} are fulfilled. Moreover, let $((\uebar,\uzbar),
 (\bmu,\bvp,\bsigma))
\in{\cal U}_R\times \Phi$ be arbitrary. Then the linear initial-boundary value problem \eqref{aux1}--\eqref{aux5}
has for every $(k_1,k_2)\in \Liq^2$ and every $(f_1,f_2,f_3)\in {L^\infty(Q)\times(H^1(0,T;H)\cap L^\infty(Q))
\times L^\infty(Q)}$ a unique solution
\,$(\mu,\phi,\sigma)\in {\cal Y}$. 
Moreover, the linear mapping $$\,((k_1,k_2),(f_1,f_2,f_3),(\mu_0,\phi_0,\sigma_0))
\mapsto (\mu,\phi,\sigma)\,$$ is continuous from {$\,\Liq^2\times \bigl(L^\infty(Q)\times (H^1(0,T;H)\cap \Liq)
\times\Liq\bigr)\times {W_{0}^3}$}
into $\,\Y$.
\Elem                 
\Bdim
We use a standard Faedo--Galerkin approximation. To this end, let $\{\lambda_k\}_{k\in\enne}\,$ and 
$\,\{e_k\}_{k\in\enne}\,$ denote the eigenvalues and associated eigenfunctions of the eigenvalue problem
$$-\Delta y+y=\lambda y\quad\mbox{in }\,\Omega,\quad \dn y=0\quad\mbox{on }\,\Gamma,$$
where the latter are normalized by $\,\|e_k\|_2=1$. Then $\,\{e_k\}_{k\in\enne}\,$ forms a complete orthonormal   
system in $\,H\,$ which is also dense in $\,V$. We put $\,V_n:={\rm span}\,\{e_1,\ldots,e_n\}$, $n\in\enne$, noting
that $\,\bigcup_{n\in\enne}V_n\,$ is dense in $\,V$. We look for functions
of the form
$$\mu_n(x,t)=\sum_{k=1}^n u_k^{(n)}(t)e_k(x), \quad \phi_n(x,t)=\sum_{k=1}^n v_k^{(n)}(t)e_k(x), \quad 
\sigma_n(x,t)=\sum_{k=1}^n w_k^{(n)}(t)e_k(x),$$
that satisfy the system
\begin{align}
\label{gs1}
&(\alpha\dt\mu_n(t),v)+(\dt\phi_n(t),v)+(\nabla\mu_n(t),\nabla v)\,=\,(z_{n1}(t),v)\quad \forall\,v\in V_n, \mbox{ for a.e. }
\,t\in (0,T),\\[0.5mm]
\label{gs2}
&(\beta\dt\vp_n(t),v)+(\nabla\vp_n(t),\nabla v)-(\mu_n(t),v)\,=\,(z_{n2}(t),v)\quad\forall\,v\in V_n, \mbox{ for a.e. }\,
\,t\in (0,T),\\[0.5mm]
\label{gs3}
&(\dt\sigma_n(t),v)+(\nabla\sigma_n(t),\nabla v)-\chi(\nabla\phi_n(t),\nabla v)\,=\,(z_{n3}(t),v)\quad\forall\,v\in V_n,
\mbox{ for a.e. }\,t\in (0,T),\\[0.5mm]
\label{gs4}
&\mu_n(0)=\lambda_4 P_n \mu_0, \quad \vp_n(0)=\lambda_4 P_n \phi_0,\quad \sigma_n(0)=\lambda_4 P_n\sigma_0,
\end{align} 
where $\,P_n\,$ denotes the $H^1(\Omega)$--orthogonal projection onto $\,V_n$, and where
\begin{align}
\label{zn1}
&z_{n1}\,=\, \lambda_1\left[P\sigma_n h(\bvp)+(P\bsigma-A-\uebar)h'(\bvp)\phi_n\right]-\lambda_2 k_1 h(\bvp)
+\lambda_3 f_1,\\
\label{zn2}
&z_{n2}\,=\,\lambda_1\left[\chi\,\sigma_n-F''(\bvp)\phi_n\right]+\lambda_3 f_2,\\
\label{zn3}
&z_{n3}\,=\,\lambda_1\left[-B\sigma_n-{E}\sigma_n h(\bvp)-
{E}\bsigma h'(\bvp)\phi_n\right]+\lambda_2 k_2+\lambda_3 f_3\,.
\end{align}
Insertion of $v=e_k$, for $k\in\enne$, in
\eqref{gs1}--\eqref{gs3},
and substitution for the second summand in \eqref{gs1} by means of \eqref{gs2}, then lead to an initial value
problem for an explicit  linear system of ordinary differential equations for the unknowns $\,u_1^{(n)},\ldots,u_n^{(n)},v_1^{(n)},\ldots,v_n^{(n)},w_1^{(n)},\linebreak$
$\ldots,w_n^{(n)}$, in which all of the coefficient
functions belong to $L^\infty(0,T)$. Hence, by virtue of Carath\'eodory's theorem, there exists a unique solution
in $W^{1,\infty}(0,T;\erre^{3n})$ that specifies the unique solution $(\mu_n,\phi_n,\sigma_n)\in W^{1,\infty}(0,T;
\Hdue)^3$ to the system \eqref{gs1}--\eqref{gs4}, for $n\in\enne$.

We now derive some a priori estimates for the Faedo--Galerkin approximations. In this procedure, $C_i>0$,
$i\in\enne$, will denote constants that are independent of $n\in\enne$ and the data $((f_1,f_2,f_3),(\mu_0,\phi_0,
\sigma_0))$, while the constant $M>0$ is given by 
\begin{align}
\label{defM}
M:=&\lambda_2\,\|(k_1,k_2)\|_{L^\infty(Q)^2}\,+\,\lambda_3\,{\|(f_1,f_2,f_3)\|_{L^\infty(Q)\times(H^1(0,T;H)\cap \Liq)\times\Liq}}\nonumber\\[0.5mm]
&+\,\lambda_4\,\|(\mu_0,\phi_0,\sigma_0)\|_{\Hdue^3}.
\end{align}
Moreover, $\,(\bmu,\bvp,\bsigma)\in \Phi$, and thus it follows that \neuj{$\bsigma, h(\bvp), h'(\bvp),F''(\bvp)\in
C^0(\overline Q)$}. Hence, there is some constant $C_1>0$ such that, for a.e. $(x,t) \in  Q$ and for all $n\in\enne$,
\begin{align}
\big(|z_{n1}|+|z_{n2}|+|z_{n3}|\big)(x,t)\,&\le\,C_1\,\big(\lambda_1(|\vp_n|+|\sigma_n|)(x,t) +
 \lambda_2(|k_1|+|k_2|)(x,t)\\
& \qquad\quad\,+ \lambda_3(|f_1|+|f_2|+|f_3|)(x,t)\big)\nonumber\\[0.5mm]
 &\,\,\le\,C_1\,\bigl(\lambda_1(|\vp_n|+|\sigma_n|)(x,t)\,+\,M\bigr).
\label{basic}
\end{align}

\vspace{1mm}\noindent
{\sc First estimate.}
\quad We insert $v=\mu_n(t)$ in \eqref{gs1}, $v=\dt\phi_n(t)$ in \eqref{gs2}, and $v=\sigma_n(t)$ in \eqref{gs3},
and add the resulting equations, whence a cancellation of two terms occurs. Then, in order to recover the
full $H^1(\Omega)-$norm below, we add to both sides of the resulting equation the same term $\,{\frac 12\frac d{dt}\|\phi_n(t)\|_2^2=}(\phi_n(t),\dt\phi_n(t))$. Integration over $[0,\tau]$, where $\tau\in (0,T]$, then
yields the identity
\begin{align}
\label{John}
&\frac 1 2\,\bigl(\alpha\|\mu_n(\tau)\|_2^2\,+\,\|\phi_n(\tau)\|_V^2 {\,+\,\|\sigma_n(\tau)\|_2^2\bigr)}\,+{\int_0^\tau\!\!\!\iO\bigl(|\nabla\mu_n|^2+|\nabla\sigma_n|^2\bigr)}
\,+\,\beta\int_0^\tau\!\!\!\iO|\dt\phi_n|^2\nonumber\\
&={\,\frac{\lambda_4^2} 2\left(\alpha\|P_n\mu_0\|_2^2+\|P_n\phi_0\|_V^2+\|P_n\sigma_0\|_2^2\right)}
\,+\int_0^\tau (\mu_n(t),z_{n1}(t))\,dt \,+{\int_0^\tau(\sigma_n(t),z_{n3}(t))\,dt}\nonumber\\
&\qquad{+\int_0^\tau(\dt\phi_n(t),z_{n2}(t)+\phi_n(t))\,dt\,+\,\chi\int_0^\tau(\nabla\phi_n(t),\nabla\sigma_n(t))\,dt
\,=:\,\sum_{i=1}^5J_i,}
\end{align}                                
{with obvious notation. We estimate the terms on the right-hand side individually. First observe that
$\|y\|_V\le\|y\|_{\Hdue}$ for all $y\in\Hdue$, and thus, 
for all $n\in\enne$,
\begin{equation}
\label{inibound}
{|J_1|\,\le\,C_2\,\lambda_4^2
\,\|(P_n\mu_0,P_n\phi_0,P_n\sigma_0)\|_{V\times V\times V}^2\,\le\,C_2\,\lambda_4^2\,
\|(\mu_0,\phi_0,\sigma_0)\|_{V\times
V\times V}^2\,\le\,C_2\,M^2.}
\end{equation}
Moreover, by virtue of \eqref{basic} and Young's inequality,
\begin{align}
\label{Paul}
{|J_2|+|J_3|\,\le\,C_3\,M^2\,+\,C_4\int_0^\tau\!\!\!\iO\left(|\mu_n|^2+|\phi_n|^2+|\sigma_n|^2\right).}
\end{align}
Likewise,
\begin{align}
\label{George}
|J_4|\,\le\,&{\,\frac{\beta}2\int_0^\tau\!\!\!\iO|\dt\phi_n|^2\,+\,\frac{C_5}{\beta}\,M^2\,+\,\frac{C_6}{\beta}
\int_0^\tau\!\!\!\iO\left(|\mu_n|^2+|\phi_n|^2+|\sigma_n|^2\right).}	
\end{align}
Finally, we have that
\begin{equation}
\label{Ringo}
{|J_5|\,\le\,\frac 12\int_0^\tau\!\!\!\iO|\nabla\sigma_n|^2\,+\,\frac {\chi^2}2\int_0^\tau\!\!\!\iO|\nabla\phi_n|^2.}
\end{equation}
Combining the estimates \eqref{John}--\eqref{Ringo}, where we subtract the first integral in \eqref{George} 
from the associated term on the left-hand side of \eqref{John}}, we have shown that
\begin{align*}
&{\frac 12\left(\alpha\|\mu_n(\tau)\|_2^2+\|\phi_n(\tau)\|_V^2+\|\sigma_n(\tau)\|_2^2\right)\,+\,\int_0^\tau\!\!\!\iO
\bigl(|\nabla\mu_n|^2\,+\,\frac 12\,|\nabla\sigma_n|^2\bigr)\,+\,\frac{\beta}2\int_0^\tau\!\!\!\iO|\dt\phi_n|^2
}\nonumber\\
&\quad{\le\,C_7\,M^2\,+\,C_8\int_0^\tau\bigl(\|\mu_n(t)\|_2^2
+\|\phi_n(t)\|_V^2+\|\sigma_n(t)\|_2^2\bigr)\,dt\,.}
\end{align*}
{Therefore, invoking
Gronwall's lemma, we conclude that, for all $n\in\enne$,}
\begin{align}
\label{esti1}
\|\mu_n\|_{L^\infty(0,T;H)\cap L^2(0,T;V)}\,+\,\|\phi_n\|_{H^1(0,T;H)\cap L^\infty(0,T;V)}
{\,+\,\|\sigma_n\|_{L^\infty(0,T;H)\cap L^2(0,T;V)}}
\,\le\,C_9 M.                                             
\end{align}  
\vspace{1mm}\noindent
{\sc Second estimate.}\quad
Next, we insert $v=\dt\mu_n(t)$ in \eqref{gs1} and integrate 
over $[0,\tau]$, where $\tau\in (0,T]$, to obtain the identity
\begin{align*}
&{\frac 12\,\|\nabla\mu_n(\tau)\|_2^2\,+\,\alpha\int_0^\tau\|\dt\mu_n(t)\|_2^2\,dt}\\
&{=\,\frac {\lambda_4^2}2\,\|\nabla P_n\mu_0\|_2^2 +\int_0^\tau(\dt\mu_n(t),z_{n1}
(t))\,dt
\,-\int_0^\tau(\dt\mu_n(t),\dt\phi_n(t))\,dt.}
\end{align*}
{Applying Young's inequality appropriately, where we make use of 
\eqref{basic} and \eqref{esti1}, we conclude the estimate}
\begin{align}\label{esti2}
&\|\mu_n\|_{H^1(0,T;H)\cap L^\infty(0,T;V)}\,\le\,C_{10}\,M.
\end{align}

\vspace{1mm}\noindent
{\sc Third estimate.} \quad At this point, we insert $v=-\Delta\mu_n(t)$ in \eqref{gs1} and $v=-\Delta\phi_n(t)$
in \eqref{gs2}, add, and integrate over $[0,\tau]$ where $\tau\in (0,T]$. We then obtain that
\begin{align*}
&\frac \alpha 2\,\|\nabla\mu_n(\tau)\|_2^2\,+\,\frac{\beta} {2}\,\|\nabla\phi_n({\tau})\|_2^2\,+\int_0^\tau\|\Delta\mu_n(t)\|_2^2\,dt
\,+\int_0^\tau\|\Delta\phi_n(t)\|_2^2\,dt \\
&=\,\frac{\alpha\lambda_4^2} 2\|\nabla P_n\mu_0\|_2^2\,+\,\frac {\beta\lambda_4^2}2 \|\nabla P_n\phi_0\|_2^2\,-\int_0^\tau(\Delta\mu_n(t),z_{n1}(t))\,dt\\
&\quad \,\,-\int_0^\tau(\Delta\phi_n(t),\mu_n(t)+z_{n2}(t))\,dt\,,
\end{align*}
whence, using \eqref{basic}--\eqref{esti2} and Young's inequality, 
\begin{equation}
\label{Hugo}
\int_0^T\bigl(\|\Delta\mu_n(t)\|_2^2\,+\,\|\Delta\phi_n(t)\|_2^2\bigr)\,dt\,\le\,C_{11}M^2\quad\forall\,n\in\enne.
\end{equation}
At this point, we invoke a classical elliptic estimate (see, e.g., \cite[Chap.~2,~Thm.~5.1]{LM}): there is a 
constant $C_\Omega>0$, which only depends on $\Omega$, such that, for every $v\in\Hdue$,
\begin{equation}\label{elli}
\|v\|_{\Hdue}\,\le\,C_\Omega\bigl(\|\Delta v\|_{\Ldue}\,+\,\|v\|_{\Huno}\,+\,\|\dn v\|_{H^{1/2}(\Gamma)}\bigr)\,.
\end{equation} 
In view of the zero Neumann boundary condition satisfied by $\mu_n$ and $\phi_n$, we thus conclude from
\eqref{esti1}, \eqref{esti2}, and \eqref{Hugo}, that
\begin{equation}\label{esti3}
\|\mu_n\|_{L^2(0,T;\Hdue)}\,+\,\|\phi_n\|_{L^2(0,T;\Hdue)}\,\le\,C_{12} M \quad\forall\,n\in\enne.
\end{equation}
With the estimate \eqref{esti3} at hand, we may (by first taking $v=\dt\sigma_n(t)$ in \eqref{gs3} and then $v=-\Delta\sigma_n(t)$) infer by similar reasoning that also
\begin{equation}\label{esti4}
\|\sigma_n\|_{H^1(0,T;H)\cap L^\infty(0,T;V)\cap L^2(0,T;\Hdue)}\,\le\, C_{13} M \quad\forall\,n\in\enne.
\end{equation}

At this point, we can conclude from standard weak and weak-star compactness arguments the existence of a triple
$(\mu,\phi,\sigma)$ such that, possibly only on  a subsequence which is again indexed by $\,n$, 
\begin{align*}
&\mu_n\to\mu,\quad \phi_n\to\phi,\quad \sigma_n\to\sigma,\\[0.5mm]
&\mbox{all weakly-star in }\,H^1(0,T;H)\cap L^\infty(0,T;V)\cap L^2(0,T;\Hdue).
\end{align*}

Standard arguments, which need no repetition here, then show that $(\mu,\phi,\sigma)$ is a strong solution to the system \eqref{aux1}--\eqref{aux5}. 
Moreover, recalling \eqref{esti1}--\eqref{esti4}, and invoking the weak sequential lower
semicontinuity of norms, we conclude that there is some $C_{13}>0$ such that
\begin{align}
\label{esti5}
\|(\mu,\phi,\sigma)\|_{(H^1(0,T;H)\cap L^\infty(0,T;V)\cap L^2(0,T;\Hdue))^3}\,\le\,C_{13}M.
\end{align}

Next, we claim that $(\mu,\phi,\sigma)\in C^0(\overline Q)^3$ and that, with a suitable $C_{14}>0$, 
\begin{align}
\label{esti6}
\|(\mu,\phi,\sigma)\|_{C^0(\overline Q)^3}\,\le\,C_{14}M. 
\end{align} 
It is easy to argue for the solution component $\,\phi$. Indeed, we have (cf. \eqref{aux2}) 
\[
\beta\dt\phi-\Delta\phi\,=\,\mu+\lambda_1(\chi\sigma-F''(\overline\phi)\phi)+\lambda_3 f_2 =:g.
\] 

\neuj{At this point, we claim that the chain rule formulas}
\neuj{
\begin{equation}\label{chain}
\dt(F''(\bvp))=F'''(\bvp)\dt\bvp,\quad\partial_{x_i}(F''(\bvp))=F'''(\bvp)\partial_{x_i}\bvp \quad\mbox{for } i=1,2,3, \quad\mbox{a.e. in }\, Q,
\end{equation}}
\neuj{are valid. Indeed, we can argue as follows: it is known from \eqref{ssbound2} that }
\neuj{
$$
r_-<r_*\le \overline\phi\le r^*<r_+ \quad\mbox{in }\,\overline Q,
$$}
\neuj{and, by {\bf (F2)} and {\bf (F3)}, that $F''\in C^3(r_-,r_+)$. Now extend $F''$ to a function $G$ defined on the whole 
of $\erre$ by putting $G(r)=F''(r)$ on  $[r_*,r^*]$, as well as $G(r)=F''(r_*)$ for $r<r_*$ \,and\, $G(r)=F''(r^*)$ for $r>r^*$.
Then $G$ is piecewise continuously differentiable on $\erre$ with $G'\in L^\infty(\erre)$. We thus may use the known chain rule
for generalized derivatives (see, e.g., \cite[Thm.~7.8]{GT}) to see that $\,\dt(G(\overline\phi))=G'(\overline\phi)\dt\overline\phi$
a.e. in $Q$. Since $G(\overline\phi)=F''(\overline\phi)$ by construction, the claim follows for the generalized time derivative, and
similar reasoning yields the assertion for the generalized space derivatives.}

\neuj{It follows from the above considerations and estimates 
that  $\,\phi\,$ solves a linear parabolic equation with zero Neumann boundary condition and 
with a right-hand side $\,g\,$ whose first two summands are known to be bounded in $(L^\infty(0,T;H)\cap L^2(0,T;V))\subset L^{10/3}(Q)$
by an expression of the form $C_{15}M$, while this holds true for the third summand in the space $L^\infty(Q)$. 
Hence, $\,g\,$ is bounded in $L^{10/3}(Q)$. Now observe that $\frac{10}3>\frac 52$ and $\phi_0\in W_0$. Therefore, we may
invoke the classical results from  \cite[Lem.~7.12]{Fredibuch} to conclude the validity of the claim for $\,\phi$.}  

For the other solution components $\,\mu\,$ and $\,\sigma$, a similar argument is not yet possible, 
since the expressions $\,\dt\phi\,$ and $\,\Delta\phi\,$ occurring in \eqref{aux1} and \eqref{aux3}, respectively,
are so far merely known to be bounded in $L^2(Q)$. In order to prove the claim by the above argument
also for $\,\mu\,$ and $\sigma$, we are now going 
to {derive additional 
bounds for $\,\dt\phi\,$ and $\,\Delta\phi$.}

To this end, \neuj{we infer from the chain rule \eqref{chain} the identity}  
\begin{equation}\label{opa1}
 \dt g = \dt\mu +\lambda_1\bigl(\chi\dt\sigma-F''(\overline\phi)\dt\phi-F'''(\overline\phi)
\dt\overline\phi\,\phi\bigr) +\lambda_3\dt f_2,
\end{equation}
where, owing to H\"older's inequality, \eqref{ssbound3}, and the continuity of the embedding $V\subset L^4(\Omega)$,
\begin{align*}
\int_0^T\!\!\!\iO|F'''(\overline\phi)
\dt\overline\phi\,\phi|^2\,\le\,K_2^2\int_0^T\|\dt\overline\phi(t)\|_4^2\,\|\phi(t)\|_4^2\,dt
\,\le\,K_2^2\,\|\dt\overline\phi\|_{L^2(0,T;V)}^2\,\|\phi\|_{L^\infty(0,T;V)}^2.
\end{align*}
Therefore, invoking \eqref{ssbound1} and \eqref{esti5},
\begin{equation}\label{opa2}
\|\dt g\|_{L^2(Q)}\,\le\,C_{16}M.
\end{equation}

At this point, we consider the linear parabolic initial-boundary value problem
\begin{align}
\label{z1}
&{\beta\dt z-\Delta z=\dt g\quad\mbox{in }\,Q,}\\[0.5mm]
\label{z2}
&{\dn z=0 \quad\mbox{on }\,\Sigma,}\\[0.5mm]
\label{z3}
&{z(0)=\beta^{-1}(\Delta\phi(0)+g(0)) \quad\mbox{in }\,\Omega,}
\end{align}
where, owing to \eqref{aux5}, 
$$\Delta\phi(0)+g(0)=\lambda_4\bigl(\Delta\phi_0+\mu_0+\lambda_1(\chi\sigma_0-F''(\phi_0)\phi_0)\bigr)+\lambda_3f_2(0)\in\Ldue.
$$
 Since also $\dt g\in L^2(Q)$, 
it follows from a classical argument that the above system admits a unique weak solution 
$\,z\in H^1(0,T;V^*)\cap C([0,T];H)\cap L^2(0,T;V)$, and since $\dn\phi_0=0$, it is easily checked that
the function
\[
w(x,t):=\phi_0(x)+\int_0^tz(x,s)\,ds \quad\mbox{for a.e. }\,(x,t)\in Q
\]
coincides with $\,\phi$, that is, in particular, we have $z=\dt\phi$. Moreover, standard estimates 
and \eqref{esti5}, \eqref{opa2} show that
\begin{align}\label{opa3}
\|\dt\phi\|_{H^1(0,T;V^*)\cap C([0,T];H)\cap L^2(0,T;V)}\,\le\,C_{17}\bigl(\|\Delta\phi(0)+g(0)\|_2\,+\,
\|\dt g\|_{L^2(Q)}\bigr)\,\le\,C_{18}M.
\end{align}
\neuj{Comparison in \eqref{aux2} and the elliptic estimate \eqref{elli} then show that also
}
\neuj{
\begin{equation}\label{oma1}
\|\phi\|_{L^\infty(0,T;W_0)}\,\le\,C_{19}M.
\end{equation}
} 

\neuj{Next, we consider the parabolic problem for $\mu$ which results if we bring the term $\dt\phi$ to the 
right-hand side of \eqref{aux1}. Arguing as above, we then find that the $L^{10/3}(Q)-$norm of the resulting
right-hand side is bounded by an expression of the form $\,C_{20}M$. Hence, owing to \cite[Lem.~7.12]{Fredibuch} again,
we can infer that $\mu\in C^0(\overline Q)$ and}
\neuj{
\begin{equation}\label{esti13}
\|\mu\|_{C^0(\overline Q) }  \,\le\,C_{21}M.
\end{equation}}

\neuj{At this point, we write \eqref{aux3} in the form}
\neuj{
$$ \dt(\sigma-\chi\phi)-\Delta(\sigma-\chi\phi)= \lambda_1[-B\sigma-E\sigma h(\overline\phi)-E\overline\sigma 
h'(\overline\phi)\phi]+\lambda_2k_2+\lambda_3f_3-\chi\dt\phi.
$$}                                           
\neuj{Apparently, the $\,L^{10/3}(Q)-$norm of the right-hand side of this equation is again
bounded by an expression of the form $C_{22}M$. Hence, we may employ \cite[Lem.~7.12]{Fredibuch} once more
to infer that $\sigma-\chi\phi\in C^0(\overline Q)$ so that also $\sigma\in C^0(\overline Q)$. In addition, we 
obtain the estimate}
\neuj{
\begin{equation}
\|\sigma\|_{C^0(\overline Q)}\,\le\,C_{23}M.
\end{equation}
}
\neuj{The above claim and \eqref{esti6} are thus shown. Finally, by recalling \eqref{esti5}, \eqref{opa3}, and \eqref{oma1},  and using 
the continuity of the embedding of $H^1(0,T;H)\cap L^2(0,T;\Hdue)$
in $C^0([0,T];V)$, we obtain that }
$$
\neuj{\|(\mu,\varphi,\sigma)\|_{\cal X}\,\le\,C_{24}M.}
$$
It then immediately follows that $(\mu,\phi,\sigma)\in\Y$, as well as
\begin{equation*}
\|(\mu,\phi,\sigma)\|_{\Y}\,\le\,C_{25}M.
\end{equation*}

With this, the existence of a solution with the asserted properties is shown. It remains to
prove the uniqueness. To this end, let $(\mu_i,\phi_i,\sigma_i)\in\Y$, $i=1,2$, be two solutions to 
the system. Then $(\mu,\phi,\sigma):=(\mu_1,\phi_1,\sigma_1)-
(\mu_2,\phi_2,\sigma_2)$ solves the system \eqref{aux1}--\eqref{aux5} 
with zero initial data, where the terms $\lambda_2k_i$,
$i=1,2$, and 
$\lambda_3f_i$, $i=1,2,3$, on the  right-hand sides do not occur. {By the definition of $\,{\cal Y}$
(recall \eqref{calX} and \eqref{calY}), and since} $(\mu,\phi,\sigma)\in\Y$,
all of the generalized partial derivatives occurring in \eqref{aux1}--\eqref{aux3} belong to $L^2(Q)$.
Therefore, we may repeat -- now for the continuous
problem -- the a priori estimates performed for the Faedo--Galerkin approximations that led us to the
estimate \eqref{esti1}. We then find analogous estimates for
$(\mu,\phi,\sigma)$, where this time the constant $\,M\,$ from \eqref{defM} equals zero. Thus, $(\mu,\phi,\sigma)
=(0,0,0)$. With this, the uniqueness is shown, which finishes the proof of the assertion.
\Edim

%\vspace{2mm}
%\Brem
%As it follows from the above proof, the solution component $\,\phi\,$ enjoys the additional regularity
%\juerg{$\,\phi\in H^1(0,T;V)$.}
%\Erem

\vspace{2mm}
With Lemma 2.4 shown, we are in a position to prepare for the application of the implicit function
theorem.  For this purpose, let us consider two auxiliary linear initial-boundary value problems. The first,
\begin{align}
\label{sysG11} 
&\alpha\dt\mu+\dt\phi-\Delta\mu\,= \, f_1\quad\mbox{in }\,Q,\\[0.5mm]
\label{sysG12} 
&\beta\dt\phi-\Delta\phi-\mu\,= \, f_2 \quad\mbox{ in \,$Q$},
\\[0.5mm]
\label{sysG13}
&\dt\sigma-\Delta\sigma+\chi\Delta\phi\,=\,f_3\quad\mbox{ in \,$Q$},\\[0.5mm]
\label{sysG14}
&\dn\mu\,=\,\dn\phi\,=\,\dn\sigma\,=\,0 \quad\mbox{ on \,$\Sigma$},\\[0.5mm]
\label{sysG15}
&\mu(0)\,=\,\phi(0)\,=\,\sigma(0)\,=0,\,\,\,\mbox{ in }\,\Omega,
\end{align}
is obtained from \eqref{aux1}--\eqref{aux5} for $\lambda_1=\lambda_2=\lambda_4=0, \, \lambda_3=1$.
Thanks to Lemma \ref{L2.4}, this system has for each $(f_1,f_2,f_3) \in L^\infty(Q)\times\bigl(
H^1(0,T;H)\cap\Liq\bigr)\times\Liq$ a unique solution
$(\mu,\phi,\sigma) \in {\cal Y}$, and the associated linear mapping
\begin{equation}\label{defG1}
{\cal G}_1:{\left(L^\infty(Q)\times\bigl(H^1(0,T;H)\cap\Liq\bigr)\times\Liq\right)}\to {\cal Y}; \,\,(f_1,f_2,f_3)\mapsto (\mu,\phi,\sigma),
\end{equation}
is continuous. The second system reads
\begin{align}\label{sysG21}
&\alpha\dt\mu+\dt\phi-\Delta\mu\,=\,0 \quad\mbox{in }\,Q,\\[0.5mm]
\label{sysG22}
&\beta\dt\phi-\Delta\phi-\mu\,=\,0 \quad\mbox{ in \,$Q$},
\\[0.5mm]
\label{sysG23}
&\dt\sigma-\Delta\sigma+\chi\Delta\phi\,=\,0\quad\mbox{ in \,$Q$},\\[0.5mm]
\label{sysG24}
&\dn\mu\,=\,\dn\phi\,=\,\dn\sigma\,=\,0 \quad\mbox{ on \,$\Sigma$},\\[0.5mm]
\label{sysG25}
&\mu(0)\,=\,\mu_0,\quad \phi(0)\,=\,\phi_0, \quad \sigma(0)\,=\,
\sigma_0,\,\,\,\mbox{ in }\,\Omega.
\end{align}
For each $(\mu_0,\phi_0,\sigma_0) \in W_0^3$, it also enjoys a unique solution $(\mu,\phi,\sigma) \in {\cal Y}$, and the associated mapping
\begin{equation}\label{defG2}
{\cal G}_2: W_0^3\to {\cal Y}; \,\,(\mu_0,\phi_0,\sigma_0)
\mapsto (\mu,\phi,\sigma),
\end{equation}
is linear and continuous as well.

\noindent In addition, we define on the open set ${\cal A}:=({\cal U}_R\times\Phi)
\subset ({\cal U}\times\Y)$ the nonlinear mapping
\begin{align}
\label{defG3}
&{\cal G}_3:{\cal U}_R\times\Phi\to {\left(L^\infty(Q)\times\bigl(H^1(0,T;H)\cap\Liq\bigr)\times\Liq\right)};
\nonumber\\[0.5mm]
&((u_1,u_2),(\mu,\phi,\sigma))\mapsto (f_1,f_2,f_3),
\,\,\mbox{ where}\nonumber\\[0.5mm]
&(f_1,f_2,f_3)=((P\sigma-A-u_1)h(\phi), \chi\,\sigma-F'(\phi),B(\sigma_s-\sigma)-{E}\sigma h(\phi)+u_2)\,. 
\end{align}
The solution $(\mu,\phi,\sigma)$ to the nonlinear state equation \eqref{ss1}--\eqref{ss5} is the sum
of the solution to the system \eqref{sysG11}--\eqref{sysG15}, where $(f_1,f_2,f_3)$ is chosen as above (with $(\mu,\varphi,\sigma)$ considered as known), and of the solution to the system
\eqref{sysG21}--\eqref{sysG25}. Therefore, the state vector $(\mu,\phi,\sigma)$ associated with the control vector $(u_1,u_2)$ is the unique solution to the nonlinear equation
\begin{equation} \label{nonlineq}
(\mu,\phi,\sigma) = {\cal G}_1 \big({\cal G}_3((u_1,u_2),(\mu,\phi,\sigma)\big) +  {\cal G}_2(\mu_0,\varphi_0,\sigma_0).
\end{equation}

Let us now define  the nonlinear mapping  $\,{\cal F}:{\cal A}\to \Y$,
\begin{align}
\label{defF}
 {\cal F}((u_1,u_2),(\mu,\phi,\sigma))\,=\,{\cal G}_1({\cal G}_3
((u_1,u_2),(\mu,\phi,\sigma))+{\cal G}_2(\mu_0,\phi_0,\sigma_0) - (\mu,\phi,\sigma).
\end{align} 
With $ {\cal F}$, the state equation can be shortly written as
\begin{equation} \label{nonlineq2}
{\cal F}((u_1,u_2),(\mu,\phi, \sigma))=(0,0,0).
\end{equation}
This equation just means that $(\mu,\phi,\sigma)$ is a solution to the state system \eqref{ss1}--\eqref{ss5} such that    $((u_1,u_2),(\mu,\phi,\sigma))\in{\cal A}$. From Theorem 2.1 we 
know that such a solution exists for every $(u_1,u_2)\in{\cal U}_R$. A fortiori, any such solution automatically enjoys the separation property  \eqref{ssbound2} and is uniquely determined.  

We are going to apply the implicit function theorem to the equation  \eqref{nonlineq2}. To this end,
we need the differentiability of the involved mappings.

Observe that, owing to the differentiability properties of the
involved Nemytskii operators (see, e.g., \cite[Thm.~4.22]{Fredibuch}), the mapping $\,{\cal G}_3\,$ is twice continuously Fr\'echet differentiable 
{as a mapping into the space $L^\infty(Q)^3$}, and for the first partial derivatives at any point 
$\,((\uebar,\uzbar),(\bmu,\bvp,\bsigma))\in {\cal A}$, and for all $(u_1,u_2)\in{\cal U}$
and $(\mu,\phi,\sigma)\in \Y$, we have the identities
\begin{align}
\label{Freu}
D_{(u_1,u_2)}{\cal G}_3((\uebar,\uzbar),(\bmu,\bvp,\bsigma))(u_1,u_2)\,=&\,
(-u_1 h(\bvp),0, u_2),\\[0.5mm]
\label{Frey}
D_{(\mu,\phi,\sigma)}{\cal G}_3((\uebar,\uzbar),(\bmu,\bvp,\bsigma))(\mu,\phi, \sigma)\,=&\,
((P\bsigma-A-\uebar)h'(\bvp)\phi+P\sigma h(\bvp), \chi\sigma-F''(\bvp)\phi,\nonumber\\
&\hspace*{2mm} -B\sigma-{E}\sigma h(\bvp)-{E}\bsigma h'(\bvp)\phi).
\end{align}  

{We now claim that the second component of ${\cal G}_3$ is also twice continuously Fr\'echet differentiable
as a mapping into $H^1(0,T;H)$. Once this will be shown, it will follow that ${\cal G}_3$ is twice continuously Fr\'echet
differentiable into the space $L^\infty(Q)\times \bigl(H^1(0,T;H)\cap L^\infty(Q)\bigr)\times L^\infty(Q)$.} 

{We first show the existence of the first derivative as a mapping into $H^1(0,T;H)$. To this end, let $((\uebar,\uzbar),(\bmu,\bvp,\bsigma))
\in {\cal U}_R\times \Phi$ and an admissible increment $(\mu,\vp,\sigma)\in{\cal Y}$ be given. Since the contribution
$\chi\sigma$ to the second component of ${\cal G}_3$ is linear, we apparently only need to show that}
\begin{equation}\label{G31}
{\lim_{\|(\mu,\vp,\sigma)\|_{\cal Y}\to 0}\,\|F'(\overline\varphi+\vp)-F'(\overline\varphi)-F''(\bvp)\vp\|
_{H^1(0,T;H)}\,=\,0.   }
\end{equation}
{In the following, we denote by $C_i$, $i\in\enne$, constants which are independent of the increments.  
Now observe that Taylor's theorem with integral remainder yields that, a.e. in $Q$,}
\begin{equation}
\label{G32}
\neuj{F^{(i)}(\bvp+\vp)-F^{(i)}(\bvp)-F^{(i+1)}(\bvp)\vp\,=\,R_i(\bvp)\vp^2,
\quad\mbox{for $i=1,2$},} 
\end{equation}
{with the remainder}
\begin{equation}
\label{G33}
{R_i(\bvp)=\int_0^1(1-s)F^{(i+2)}(\bvp+s\vp)\,ds\,.}
\end{equation}
\neuj{At this point, we recall that $F\in C^5(r_-,r_+)$ and the global bounds \eqref{ssbound2} and \eqref{ssbound3}. 
We thus may argue as in the derivation of \eqref{chain}, using the chain rule \cite[Thm.~7.8]{GT} for generalized 
derivatives, to conclude that it holds the estimate}
\begin{equation}\label{rbound}
{|R_i(\bvp)|\le C_1, \quad|\dt( R_i(\bvp))|\le C_1\bigl(|\dt\bvp|+|\dt\vp|\bigr),\quad\mbox{a.e. in \,$Q$, for $i=1,2,3$.}}
\end{equation}
{We thus have }
{
\begin{align}\label{G34}
&\|R_1(\bvp)\vp^2\|_{H^1(0,T;H)}^2 \le \int_Q\bigl(|\dt (R_1(\bvp))|^2|\vp|^4+|2\,R_1(\bvp)\vp\dt\vp|^2 + |R_1(\bvp)|^2 |\vp|^4\bigr)
\nonumber\\
&\le C_2\,\|\vp\|_{L^\infty(Q)}^4\Big(1+\int_Q\bigl(|\dt\bvp|^2+|\dt\vp|^2\bigr)\Big)\,+\,C_3
\,\|\vp\|_{L^\infty(Q)}^2\int_Q|\dt\vp|^2\non\\[1mm]
&\le C_4\,\|\vp\|_{\widetilde X}^4\,\le\,C_4\,\|(\mu,\phi,\sigma)\|_{\cal Y}^4.
\end{align}}
{Hence, \eqref{G31} is shown, which proves the claim for the first Fr\'echet derivative. A similar calculation
along the same lines, which may be omitted here, then yields that also the second Fr\'echet derivative of the
second component of ${\cal G}_3$ exists and that the corresponding contribution can
be expressed by $F^{(3)}(\bvp)\phi_1\phi_2$ for the directions $(\mu_1,\phi_1,\sigma_1), (\mu_2,\phi_2,\sigma_2)\in{\cal Y}$.} 

{It remains to
show that we have continuous second-order Fr\'echet differentiability. To this end, let the increment
${\bf h}=(h_1,h_2,h_3)\in{\cal Y}$ be
such that $(\bmu+h_1,\bvp+h_2,\bsigma+h_3)\in\Phi$. Moreover, let the directions $(\mu_i,\phi_i,\sigma_i)\in{\cal Y}$ satisfy
$\|(\mu_i,\phi_i,\sigma_i)\|_{\cal Y}=1$ for $i=1,2$. Then, in particular, $\|\phi_i\|_{L^\infty(Q)}\le 1, \,\,i=1,2$.
Moreover,}
\neuj{
\begin{align*}
F^{(3)}(\bvp+h_2)-F^{(3)}(\bvp)=h_2\,R(\bvp) :=h_2\int_0^1 F^{(4)}(\bvp+sh_2)\,ds,
\end{align*}}
\neuj{where, arguing as above,}
\neuj{
$$
|R(\bvp)|\le C_5,\quad\left|\dt (R(\bvp))\right|\le C_6\left(|\dt\bvp|+|\dt\phi|\right),\quad\mbox{a.e. in }\,Q.
$$} 
Therefore, using the fact that $\|\phi_i\|_{L^\infty(Q)}\le 1, \,\,i=1,2,$ we infer that
{
\begin{align*}
&\left\|\bigl(F^{(3)}(\bvp+h_2)-F^{(3)}(\bvp)\bigr)\phi_1\phi_2\right\|^2_{H^1(0,T;H)}\\
&\le\,C_7\int_Q |h_2|^2\left(|\dt\bvp|^2+|\dt h_2|^2\right)\,+\,C_8\int_Q|\dt h_2|^2\,+\,C_9\int_Q|h_2|^2\left(1+
|\dt\phi_1|^2+|\dt\phi_2|^2\right)\\[1mm]
&\le C_{10}\left(\|h_2\|_{L^\infty(Q)}^2\,+\,\|\dt h_2\|^2_{L^2(Q)}\right)\,\le\,C_{10}\,
\|(h_1,h_2,h_3)\|^2_{\cal Y}\,.
\end{align*}}
{Hence, we even have local Lipschitz continuity. With this, the above claim is finally proved.}

\vspace{2mm}
{At this point, we may conclude from the chain rule  that ${\cal F}$} 
is twice continuously Fr\'echet differentiable from ${\cal U}_R\times\Phi$ into $\Y$, with the first-order partial derivatives
\begin{align}
D_{(u_1,u_2)}{\cal F}((\uebar,\uzbar),(\bmu,\bvp,\bsigma))\,&=\,{\cal G}_1\circ D_{(u_1,u_2)}{\cal 
 G}_3((\uebar,\uzbar),(\bmu,\bvp,\bsigma)),
\label{DFu}
\\[1mm]
D_{(\mu,\phi,\sigma)}{\cal F}((\uebar,\uzbar),(\bmu,\bvp,\bsigma))\,&=\,{\cal G}_1\circ
D_{(\mu,\phi,\sigma)}{\cal G}_3((\uebar,\uzbar),(\bmu,\bvp,\bsigma))-I_{\Y},
\label{DFy}
\end{align}
where $\,I_{\Y}\,$ denotes the identity mapping on $\,\Y$.

We now introduce for convenience abbreviating denotations, namely,
\begin{align}
&\bu:=(u_1,u_2),\quad \overline\bu:=(\uebar,\uzbar), \quad {\bf y}:=(\mu,\vp,\sigma),\quad
{\bf {\overline y}}:=(\bmu,\bvp,\bsigma),\nonumber\\
&{\bf y}_0:=(\mu_0,\phi_0,\sigma_0), \quad \mathbf{0}:=(0,0,0).
\end{align}
With these denotations, we want to prove the differentiability of the control-to-state mapping $\bu \mapsto  {\bf y}$  defined implicitly by the equation $\,{\cal F}(\bu, {\bf y})=\mathbf{0}$, using the implicit function
theorem. Now let $\overline\bu\in  {\cal U}_R$ be given and $\overline{\bf y}={\cal S}(\overline\bu)$. We need to show that the linear
and continuous operator $\,D_{{\bf y}}{\cal F}(\overline\bu,\overline{\bf y})$ is a topological isomorphism from $\Y$ into itself. 

To this end, let ${\bf v}\in\Y$ be arbitrary. Then the identity $\,D_{\bf y}{\cal F}(\overline\bu,
\overline{\bf y})({\bf y})={\bf v}\,$ just means that $\,{\cal G}_1\left(D_{\bf y}{\cal G}_3(\overline\bu,
\overline{\bf y})(\bf y)\right)-{\bf y}={\bf v}$, which is equivalent to saying that   
\begin{equation*}
{\bf w}\,:=\,{\bf y}+{\bf v}={\cal G}_1\left(D_{\bf y}{\cal G}_3(\overline\bu,
\overline{\bf y})(\bf w)\right)-{\cal G}_1\left(D_{\bf y}{\cal G}_3(\overline\bu,
\overline{\bf y})(\bf v)\right).
\end{equation*}
The latter identity means that ${\bf w}$ is a solution to the system \eqref{aux1}--\eqref{aux5} for $\lambda_1=\lambda_3=1,
\lambda_2=\lambda_4=0$, with the specification $(f_1,f_2,f_3)=-{\cal G}_1\left(D_{\bf y}{\cal G}_3(\overline\bu,
\overline{\bf y})(\bf v)\right)\in\Y$. By Lemma 2.4, such a solution ${\bf w}\in\Y$ exists and is uniquely 
determined. We thus can infer that $D_{{\bf y}}{\cal F}(\overline\bu,\overline{\bf y})$ is surjective. At the same time, taking ${\bf v}=\mathbf{0}$, we see that the equation $D_{{\bf y}}{\cal F}(\overline\bu,\overline{\bf y})({\bf y})=\mathbf{0}$
means that ${\bf y}$ is the unique solution to \eqref{aux1}--\eqref{aux5} for $\lambda_1=1,\lambda_2=\lambda_3=
\lambda_4=0$. Obviously, ${\bf y}=\mathbf{0}$, which implies that $D_{{\bf y}}{\cal F}(\overline\bu,\overline{\bf y})$ is
also injective and thus, by the open mapping principle, a topological isomorphism from $\Y$ into itself. 

At this point, we may employ the implicit function theorem (cf., e.g., \cite[Thms. 4.7.1 and 5.4.5]{cartan1967} or \cite[10.2.1]{Dieu}) to conclude    
that the mapping $\S$ is twice continuously Fr\'echet differentiable from ${\cal U}_R$ into $\Y$ and that the
first Fr\'echet derivative $\,D\S(\ubar)\,$ of $\,\S\,$ at $\overline\bu\in{\cal U}_R$ is given by the formula
\begin{equation}
\label{DS1}
D{\cal S}(\overline\bu)\,=\,-D_{\bf y}{\cal F}(\overline\bu,\overline{\bf y})^{-1}\circ D_{\bu}
{\cal F}(\overline\bu,\overline{\bf y}).
\end{equation}
Now let ${\bf k}=(k_1,k_2)\in{\cal U}$ be arbitrary and ${\bf y}=(\mu,\phi,\sigma)=D{\cal S}(\overline\bu)({\bf k})$. 
Then, $$ D_{\bf y}{\cal F}(\overline\bu,\overline{\bf y})(\bf y)=- D_{\bu}
{\cal F}(\overline\bu,\overline{\bf y})({\bf k}),$$ which is obviously equivalent to saying that
$$
{\bf y}={\cal G}_1\left(D_{\bf y}{\cal G}_3(\overline\bu,\overline{\bf y})({\bf y})\right)
+{\cal G}_1(-k_1 h(\bvp),0,k_2).
$$
This, in turn, means that ${\bf y}$ is the unique solution to the problem \eqref{aux1}--\eqref{aux5} for
$\lambda_1=\lambda_2=1,\lambda_3=\lambda_4=0$.

In summary, we have shown the following result.

\Bthm
Suppose that the conditions {\bf (F1)}--{\bf (F4)}, {\bf (A1)}--{\bf (A3)}, {\bf (A6)} and 
\eqref{defUR} are fulfilled,
let $\overline \bu=(\uebar,\uzbar)\in {\cal U}_R$ be arbitrary and $(\bmu,\bvp,\bsigma)={\cal S}(\overline\bu)$. Then the control-to-state 
operator $\S$ is twice continuously Fr\'echet differentiable at $\,\overline\bu\,$ as a mapping from $\,{\cal U}\,$ into 
$\,{\cal Y}$. Moreover, for every $(k_1,k_2)\in {\cal U}$, the Fr\'echet derivative $\,D \S(\overline\bu)\in 
{\cal L}({\cal U},\Y)\,$ of $\,\S\,$ at
$\,\overline\bu\,$ is given by the identity $\,D\S(\overline\bu)(k_1,k_2)=(\mu,\phi,\sigma)$, where 
$(\mu,\phi,\sigma)$ is the 
unique solution to the linear system \eqref{aux1}--\eqref{aux5} with $\lambda_1=\lambda_2=1,\lambda_3=\lambda_4=0$.
\Ethm

%%%%%%%%%%%%%%%%%%%%%%%%%%%%
\section{First-order necessary optimality conditions}
%%%%%%%%%%%%%%%%%%%%%%%%%%%
\setcounter{equation}{0}

{In this section, we aim at deriving  associated first-order necessary optimality conditions. To this end, 
we define the (control) reduced objective functional 
$\,\widetilde {\cal J}\,$ by
\begin{equation}\label{reduced}
\widetilde{\cal J}(\bu) = {\cal J}(\S(\bf u),{\bf u}),
\end{equation}
where we recall that $\S(\bu)=(\mu,\varphi,\sigma)$ is the unique solution to the state system associated with
$\,\bu$. The functional $\,\widetilde{\cal J}\,$ is the sum of a nonconvex functional
$\,{\cal J}_1\,$ and the convex functional $\,\kappa\,g$, namely
\[
\widetilde {\cal J} = {\cal J}_1 + \kappa g,
\]
where} 
\begin{equation}
\label{defJ1}
{{\cal J}}_1(\bu)=  \frac{\beta_1}2 \intQ |\vp_{\bu}-\widehat \vp_Q|^2
  + \frac{\beta_2}2 \iO |\vp_{\bu}(T)-\widehat\vp_\Omega|^2
  + \frac{\nu}2 \intQ |{\bf u}|^2.
\end{equation}
Here, $\,g\,$ is one of the functionals \eqref{sparse}--\eqref{dxsparse}, and 
we denote by $\vp_{\bu}$ the second component of $\S(\bu)$.

{Since, owing to \cite[Thm.~2.6]{CSS1}, the control-to-state mapping is Fr\'echet differentiable from $L^2(Q)^2$
into $L^2(Q)\times C^0([0,T];\Ldue)\times L^2(Q)$, in particular, the functional ${\cal J}_1$ is  a Fr\'echet 
differentiable mapping from $L^2(Q)^2$ into $\erre$. Therefore, the chain rule shows that, for every $\ubar
=(\overline u_1,\overline u_2)\in L^2(Q)^2$ and $\mathbf{k}=(k_1,k_2)\in L^2(Q)^2$, it holds that}
\begin{align}
\label{DJ}
{D{\cal J}_1(\ubar)(\mathbf{k})}\,&{=\,\beta_1\int_Q(\phi_{\ubar}-\widehat\phi_Q)\,\phi\,+\,\beta_2\iO
(\phi_{\ubar}(T)-\widehat\phi_\Omega)\,\phi(T)\,+\,\nu\int_Q\ubar\cdot\mathbf{k}},
\end{align}
where $\,\phi\,$ is the second component of the solution $\,(\mu,\phi,\sigma)\,$ to the linearized system 
\eqref{aux1}--\eqref{aux5} with $\lambda_1=\lambda_2=1, \lambda_3=\lambda_4=0$, and where $"\cdot"$ stands for the euclidean inner product in $\erre^2$. 

Now assume that  $\overline \bu = (\overline u_1,\overline u_2)$ is a locally optimal control for ($\mathcal{CP}$). 
Then it is easily seen that the variational inequality
\begin{equation}
\label{var1}
{D{\cal J}_1(\ubar)(\bu - \overline \bu) + \kappa (g(\bu) - g(\overline \bu)) \ge 0 \quad \forall\, \bu \in \uad}
\end{equation}
is satisfied. Indeed, if $\bu\in\uad$ and $t\in(0,1)$ are given, then we can infer from the convexity of $\,g\,$
that
\begin{align*}
{0}\,&{\le\,{\cal J}_1(\ubar+t(\bu-\ubar))+\kappa g(\ubar+t(\bu-\ubar))-{\cal J}_1(\ubar)-\kappa g(\ubar)}\\[1mm]
&{\le\,{\cal J}_1(\ubar+t(\bu-\ubar))-{\cal J}_1(\ubar)+\kappa\,t(g(\bu)-g(\ubar)),}
\end{align*}
whence, dividing by $t>0$ and then taking the limit as $t\searrow0$, \eqref{var1} follows. But \eqref{var1} 
implies that $\ubar$ solves the convex minimization problem
\[
{\min_{\bu \in L^2(Q)^2} \,\bigl( \Phi(\bu)+ \kappa g(\bu) + I_{\uad}(\bu)\bigr),}
\]
with $\Phi(\bu)=D{\cal J}_1(\overline \bu)\bu$, and where $I_{\uad}$ denotes the indicator function of $\uad$. Hence,
denoting by the symbol $\,\partial\,$ the subdifferential mapping in $L^2(Q)^2$, we have that
\[
{\mathbf{0}\in\partial\bigl(\Phi+\kappa g+I_{\uad}\bigr)(\ubar).}
\]
At this point, we anticipate that we shall see in the next section that 
$\partial g(\bu) \subset L^2(Q)^2$ for all of our 
choices of $g$. Therefore, we may infer from the well-known rules for subdifferentials of convex functionals that
\[
{ \mathbf{0}\in \{D{\cal J}_1(\ubar)\}+\kappa\partial g(\ubar)+\partial I_{\uad}(\ubar).}
\]
In other words, there
are $\overline\bl\in\partial g(\ubar)$ and $\hat\bl\in\partial I_{\uad}(\ubar)$ such that
\begin{equation}\label{supi}
{\mathbf{0}=D{\cal J}_1(\ubar)+\kappa\overline\bl+\hat\bl.}
\end{equation}
But, by the definition of $\partial I_{\uad}(\ubar)$, we have $\,\hat\bl(\bu-\ubar)\le 0$ for every
$\bu\in \uad$. Hence, thanks to \eqref{supi},
\[
{0\le D{\cal J}_1(\ubar)(\bu-\ubar)+\kappa\overline\bl(\bu-\ubar)\quad\forall\,\bu\in\uad.}
\] 
We have thus shown the following result (where we identify $\overline\bl$  
with the corresponding element of $L^2(Q)^2$ 
according to the Riesz isomorphism).
 
\Blem  If {$\ubar\in\uad$} is a locally optimal control for {\rm ($\mathcal{CP}$)}, then there is some $\overline \bl \in \partial {g(\ubar)}\subset L^2(Q)^2$ such that 
\begin{equation} \label{varineq1}
D{\cal J}_1(\overline \bu)(\bu - \overline \bu) + \kappa \int_Q \overline
\bl(x,t)\cdot (\bu(x,t)- \overline \bu(x,t))\, dx\,dt \ge 0 \quad \forall \,\bu \in \uad.
\end{equation}
\Elem

{
\Brem The idea for the proof of the above lemma goes back to \cite{ekeland_temam1974} and to the papers \cite{casas_ryll_troeltzsch2013,casas_ryll_troeltzsch2015}, where it has been worked out for control problems 
with semilinear reaction-diffusion equations.
The concrete form of $\partial g$ depends on the particular choice of $g$ and will be presented below.
\Erem
}

Next, we aim to simplify the expression $D{\cal J}_1(\overline \bu)(\bu - \overline \bu)$ in \eqref{varineq1}
by introducing an adjoint state.
To this end, we consider the following adjoint system:
\begin{align}
\label{adj1}
&-\alpha \partial_t \psi_1 - \Delta \psi_1 \,=\, \psi_2 \quad\mbox{ in }\,Q,\\[1ex]
\label{adj2}
&-\partial_t (\psi_1+\beta\psi_2) -\Delta (\psi_2 - \chi \psi_3)\,=\,\beta_1(\overline \varphi - \varphi_Q) 
+ (P\overline \sigma - A -\overline u_1)h'(\overline \varphi)\psi_1\nonumber\\[0.5mm]
&\hspace*{62mm}- F''(\overline \varphi)\psi_2 - {E} \overline \sigma h'(\overline \varphi)\psi_3\quad\mbox{ in }\,Q,\\[1ex]
\label{adj3}
&-\partial_t \psi_3 - \Delta \psi_3\,=\, P h(\overline \varphi)\psi_1 + \chi \psi_2 
- B\psi_3 - {E} h(\overline \varphi) \psi_3\quad\mbox{ in }\,Q,\\[1ex]
\label{adj4}
&\partial_n \psi_1 = \partial_n \psi_2 = \partial_n \psi_3 =0 \quad \mbox{ on }\,\Sigma,\\[1ex]
\label{adj5}
&\psi_1(T) =\psi_3(T) = 0,\quad  \beta\psi_2(T)= \beta_2(\overline \varphi(T) - \varphi_\Omega), \quad \mbox{ in }\,\Omega.
\end{align}

\noindent
{According to \cite[Thm.~2.8]{CSS1}, the adjoint system \eqref{adj1}--\eqref{adj5} has under the general assumptions
{\bf (F1)}--{\bf (F4)} and {\bf (A1)}--{\bf (A6)} a unique weak solution 
$\boldsymbol{\psi}=(\psi_1,\psi_2,\psi_3)$ with the regularity}
\begin{align}
\label{reg1}
&{\psi_1\in H^1(0,T;H)\cap C^0([0,T];V)\cap L^2(0,T;W_0),}\\
\label{reg23}
&{\psi_2,\psi_3\in H^1(0,T;V^*)\cap C^0([0,T];H)\cap L^2(0,T;V).}
\end{align}
{We have the following result.}
\Bthm
\,\,{\rm (Necessary optimality condition)} \,\,\,{Suppose that {\bf (F1)}--{\bf (F4)} and {\bf (A1)}--{\bf (A6)}
are fulfilled, and let $\overline \bu \in \uad$ be} a locally optimal control of  {\rm ($\mathcal{CP}$)} 
{with associated state $\,(\overline\mu,\overline\phi,\overline\sigma)={\cal S}(\ubar)$
and} adjoint state $\,\overline {\boldsymbol{\psi} }= (\overline \psi_1,\overline \psi_2,\overline \psi_3)$.
 {Then, there exists} some $\overline {\boldsymbol{\lambda}} = 
 (\lambda_1,\lambda_2)^\top \in \partial g(\overline \bu)$ such that 
\begin{equation}\label{varineq2}
\int_Q  \left( \overline {\mathbf d}(x,t) + \kappa \overline {\bl}(x,t) +\nu \overline \bu(x,t)\right)
\cdot (\bu(x,t) - \overline \bu(x,t))\, dxdt \ge 0 \quad \forall \bu \in \uad,
\end{equation}
where $\overline {\mathbf  d} \in {L^2(Q)^2}$ is defined by 
\[
\overline {\mathbf d}(x,t) = \left(
\begin{array}{c}
- \overline{\mathbf \psi}_1(x,t) h({\overline \varphi}(x,t))\\
\overline\psi_3(x,t)
\end{array}
\right) \quad\mbox{for a.e. $(x,t)\in Q$}.
\]
\Ethm
\Bdim
Using the adjoint state $\overline{\boldsymbol \psi}$, we obtain the representation 
\[
{D{\cal J}_1}(\overline \bu)(\bu - \overline \bu) = \int_Q \left( \overline {\mathbf d} +\nu
 \overline \bu\right)\cdot (\bu - \overline \bu)\, dx\,dt.
\]
This follows from the proof of \cite[Thm.~2.9]{CSS1},  where the notation $\boldsymbol \psi = (p,q,r)$
and $\bu = (u,w)$ is used. The claim is now an immediate {consequence} of \eqref{varineq1}.
\Edim

%%%%%%%%%%%%%%%%%%%%%%%%%%%%%%%%%%%%%%%%%%
\section{Sparsity of optimal controls}
%%%%%%%%%%%%%%%%%%%%%%%%%%%%%%%%%%%%%%%%%%
\setcounter{equation}{0}

The convex function $g$ in the objective functional accounts for the sparsity of optimal controls, i.e., 
the optimal control can vanish in some region of the space-time cylinder $Q$. The form of this region depends on the particular {choice} of the functional $g$, while its size depends on the sparsity parameter $\kappa$.
These sparsity properties can be deduced from the variational inequality \eqref{varineq2} and the particular
form of the subdifferential  $\partial g$. 

Therefore, we first provide known results on the subdifferential and apply them to the analysis of an auxiliary variational inequality.

%%%%%%%%%%%%%%%%%%%%%%%%%%%%%%%%%%%%%%%%%%
\subsection{Preliminaries}
%%%%%%%%%%%%%%%%%%%%%%%%%%%%%%%%%%%%%%%%%%

Let us begin with the subdifferential of the $L^2$-norm,
\[
\gamma(v) = \|v\|_{L^2(\Omega)} = \left(\int_\Omega |v(x)|^2 \,dx\right)^{1/2},
\]
which is given by (see, e.g., \cite{ioffe_tikhomirov1979}) 
\begin{equation}\label{dgamma}
\partial \gamma(v) = 
\left\{
\begin{array}{ll}
\{z \in\Ldue : \|z\|_\Ldue \le 1\} &\,\mbox{ if }\, v = 0\\[1ex]
v/ \|v\|_{\Ldue} &\,\mbox{ if }\, v \not= 0
\end{array}
\right. 
\end{equation}
{In order to} have directional sparsity with respect to time, we use the functional
\begin{align} 
&g_T: L^1(0,T;L^2(\Omega)) \to \mathbb{R},\nonumber\\[1mm]
\label{g1}
&g_T(u) = \int_0^T \left(\int_\Omega |u(x,t)|^2\, dx\right)^{1/2}dt = \int_0^T \gamma(u(t))\,dt.
\end{align}
The associated subdifferential is given by (cf., \cite{herzog_stadler_wachsmuth2012})
\[
\partial g_T(u) = 
\{\lambda \in L^\infty(0,T;L^2(\Omega)): \lambda(\cdot,t) \in \partial \gamma (u(\cdot,t)) \mbox{ for a.a. } t\in (0,T)\},
\]
that is,
\begin{equation}
\label{dg1}
\partial g_T(u) = \left\{\lambda \in L^\infty(0,T;L^2(\Omega)): 
\left\{
\begin{array}{ll}
\|\lambda(\cdot,t)\|_{\Ldue} \le 1 &\,\mbox{ if } u(\cdot,t) = 0\\[1ex]
\lambda(\cdot,t)\,=u(\cdot,t)/ \|u(\cdot,t)\|_{\Ldue}&\,\mbox{ if } u(\cdot,t) \not= 0
\end{array}
\right.
\right\},
\end{equation}
where the properties above are satisfied for a.a. $t \in (0,T)$.

Directional sparsity with respect to space is obtained from the functional
\begin{align} 
&g_\Omega: L^1(\Omega;L^2(0,T)) \to \mathbb{R}, \nonumber\\[1mm]
\label{g2}
&g_\Omega(u) = \int_\Omega \left(\int_0^T |u(x,t)|^2\, dt\right)^{1/2}dx = \int_\Omega \|u(x,\cdot)\|_{L^2(0,T)}\,dx\,.
\end{align}
Interchanging the roles of $t$ and $x$, we get
\begin{equation}
\label{dg2}
\partial g_\Omega(u) = \left\{\lambda \in L^\infty(\Omega;L^2(0,T)): 
\left\{
\begin{array}{ll}
\|\lambda(x,\cdot)\|_{L^2(0,T)} \le 1 &\,\mbox{ if } u(x,\cdot) = 0\\[1ex]
\lambda(x,\cdot)=u(x,\cdot)/ \|u(x,\cdot)\|_{L^2(0,T)}&\,\mbox{ if } u(x,\cdot) \not= 0
\end{array}
\right.
\right\}
\end{equation}
where the properties above have to be fulfilled for a.a. $x \in \Omega$.

In the case of full sparsity, i.e., for 
\begin{equation} \label{g3}
g_Q: L^1(Q) \to \mathbb{R}, \quad
g_Q(u) = \int_Q |u(x,t)|\, dx\,dt,
\end{equation}
the subdifferential is classical. We have (see \cite{ioffe_tikhomirov1979})
\begin{equation}\label{dgQ}
\partial g_Q(u) = \left\{\lambda \in L^\infty(Q):\,
\lambda(x,t) \in \left\{
\begin{array}{ll}
\{1\} & \mbox{ if } u(x,t) > 0\\
{[-1,1]}& \mbox{ if } u(x,t) = 0\\
\{-1\} & \mbox{ if } u(x,t) < 0\\
\end{array}
\right.
\mbox{ for a.e. } {(x,t) \in Q}
\right\}.
\end{equation}
Here, we will concentrate on directional sparsity in time, since this seems to be the most important sparsity for medical applications. In this case, if an application to medication is considered, directional sparsity will allow to 
stop the administration of drugs in certain intervals of time.
To this {end, we now discuss} the following auxiliary variational inequality:
\begin{equation} \label{varineqaux}
\int_Q (d(x,t) + \kappa \lambda(x,t) + \nu u(x,t))(v(x,t) - u(x,t))\,dx\,dt \ge 0 \quad \forall\, v \in C,
\end{equation}
where $\lambda \in \partial g_T(u)$ and 
\begin{equation} \label{C}
C = \{v \in L^\infty(Q): \underline u \le v(x,t) \le {\hat u} \mbox{ a.e. in } Q\}
\end{equation}
with given real numbers $\underline u < 0 < \hat u$, $\kappa> 0$, $\nu > 0$, and a given function
$d \in L^2(Q)$. 
{Obviously, \eqref{varineqaux} just means that $\,u\,$ is the $\,L^2(Q)-$orthogonal projection of $\,-\frac 1\nu(d+\kappa\lambda)\,$ onto the closed and convex subset $\,C\,$ of $\,L^2(Q)$, which is well known to 
be given by the formula
\begin{equation}
\label{pro1}
u(x,t) = \mathbb{P}_{[\underline u , {\hat u}]} (-{\nu}^{-1} (d(x,t) + \kappa \lambda(x,t)))
\quad\mbox{for a.e. \,$(x,t)\in Q$},
\end{equation}
where we denote by $\mathbb{P}_{[\underline u , {\hat u}]}:  \mathbb{R} \to [\underline u , {\hat u}]$ the pointwise projection function
\begin{equation}\label{pro2}
\mathbb{P}_{[\underline u , {\hat u}]}(s) = \min\{{\hat u},\max\{\underline u,s\}\}.
\end{equation}
Moreover, it is well known that the following
pointwise relations hold true {for almost all $(x,t) \in Q$}:}     
\begin{equation} \label{cases}
\begin{array}{l}
d(x,t) + \kappa \lambda(x,t) + \nu u(x,t) > 0 \quad \Longrightarrow  \quad u(x,t) = \underline{u}\\[1ex]
d(x,t) + \kappa \lambda(x,t) + \nu u(x,t) < 0  \quad \Longrightarrow  \quad u(x,t) = {\hat u}.
\end{array}
\end{equation}

The next result is already known from \cite{herzog_stadler_wachsmuth2012,casas_herzog_wachsmuth2017}.
Nevertheless, we present a proof for the readers' convenience.
\begin{lemma}\label{L4.1} \,\,{\rm (Sparsity)} \quad Let $\,u \in C\,$ be a solution {to} the variational 
inequality \eqref{varineqaux}. Then, for a.e. $t\in (0,T)$,
\begin{equation}\label{equivalence}
u(\cdot,t) = 0 \quad \Longleftrightarrow \quad  \|d(\cdot,t)\|_{L^2(\Omega)} \le \kappa,
\end{equation}
as well as 
\begin{equation} \label{subdiff}
\lambda(\cdot,t) \left\{
\begin{array}{lcl}
\in B(0,1)&\,\mbox{ if } \, \|u(\cdot,t) \|_{L^2(\Omega)} = 0\\[1ex]
= \displaystyle \frac{u(\cdot,t)}{\|u(\cdot,t) \|_{L^2(\Omega)}}&\,\mbox{ if } \, \|u(\cdot,t) \|_{L^2(\Omega)} \not= 0
\end{array}
\right. ,
\end{equation}
where $B(0,1) = \{ v \in L^2(\Omega): \|v\|_{L^2(\Omega)}\le 1\}.$
\end{lemma}
\begin{proof}
(i) \,We first show that, for a.e. $\,t \in (0,T)$, the condition \,$\|u(\cdot,t)\|_{\Ldue} = 0\,$ implies that $\, \|d(\cdot,t)\|_{L^2(\Omega)} \le \kappa$. 
So consider the set $E = \{t\in(0,T):  \|u(\cdot,t)\|_{\Ldue} = 0\}$. Then \eqref{cases} yields that 
\[d(\cdot,t) + \kappa \lambda(\cdot,t) + 0 = 0,
\]
for a.e. $t \in E$, since otherwise the set of points $x \in \Omega$, where $u(x,t) = \underline{u}$ or $u(x,t) = 
\hat u$, would have positive measure, which contradicts the assumption that $\|u(\cdot,t)\|_{\Ldue} = 0$. 

From the equation above, we deduce that $d(\cdot,t) = - \kappa \lambda(\cdot,t)$, and thus
\[
\|d(\cdot,t)\|_{L^2(\Omega)} = \kappa \|\lambda(\cdot,t)\|_{L^2(\Omega)} \le \kappa,
\]
thanks to the form of $\partial g_T(u)$. 

\noindent(ii) Next, we confirm that the reverse implication 
\[\|d(\cdot,t)\|_{L^2(\Omega)} \le \kappa \ \Longrightarrow 
\|u(\cdot,t)\|_{\Ldue} = 0\]
holds true for almost every $t \in (0,T)$. {To this end, let 
\[
E=\{t\in (0,T): \,\|d(\cdot,t)\|_{L^2(\Omega)} \le \kappa \,\mbox{ and }\,\|u(\cdot, t)\|_{\Ldue} \not=0\}.
\]
We have to show that the Lebesgue measure $\,|E|\,$ of $\,E\,$ is zero. 
We denote by $\Omega_+(t)$ and $\Omega_-(t)$ the sets of points $\,x \in \Omega\,$ where 
$u(x,t) > 0$ and $u(x,t) < 0$, respectively.
Now recall that the implications \eqref{cases} must be satisfied. Since, by assumption, $\underline u<0<\hat u$,
we readily deduce that}
\begin{align}
\label{yeah}
&d(x,t) + \kappa \lambda(x,t) + \nu u(x,t) \le 0 \quad \mbox{for a.e. } x \in \Omega_+(t),\nonumber\\[0.5mm]
&d(x,t) + \kappa \lambda(x,t) + \nu u(x,t) \ge 0 \quad \mbox{for a.e. } x \in \Omega_-(t).
\end{align}
In $E$, we have $\|u(\cdot,t)\|_{L^2(\Omega)}\not=0$, and therefore, by \eqref{dgamma},   $\lambda(\cdot,t)= u(\cdot,t)/\|u(\cdot,t)\|_{L^2(\Omega)}$. Now the upper inequality in \eqref{yeah}
implies that
\[
d(x,t) \le - \kappa \,\frac{u(x,t)}{\|u(\cdot,t)\|_{L^2(\Omega)}} - \nu u(x,t) \quad \mbox{for a.e. } x \in \Omega_+(t). 
\]
Since both summands on the right-hand side are negative, we have
\[
{|d(x,t)| >  \kappa \,\frac{u(x,t)}{\|u(\cdot,t)\|_{L^2(\Omega)}}} \quad \mbox{for a.e. } x \in \Omega_+(t).
\]
In the same way, we deduce from the lower inequality in \eqref{yeah} that
\[
d(x,t) \ge - \kappa \frac{u(x,t)}{\|u(\cdot,t)\|_{L^2(\Omega)}} - \nu u(x,t) \quad \mbox{for a.e. } x \in \Omega_-(t),
\]
where both summands on the right-hand side are positive. This, in turn, yields that
\[
{|d(x,t)| >  \kappa \,\frac{|u(x,t)|}{\|u(\cdot,t)\|_{L^2(\Omega)}}} \quad \mbox{for a.e. } x \in \Omega_-(t).
\]
{Since $\,u(\cdot,t)\,$ vanishes on $\,\Omega\setminus(\Omega_+(t)\cup\Omega_-(t))$, we thus can infer that}
\[
\begin{aligned}
\|d(\cdot,t)\|_{L^2(\Omega)} &\ge \left(\int_{\Omega_+(t)\cup \Omega_-(t)} |d(x,t)|^2\,dx\right)^{\frac{1}{2}}> 
\kappa  \left(\int_{\Omega_+(t)\cup \Omega_-(t)} \frac{|u(x,t)|^2}{\|u(\cdot,t)\|_{L^2(\Omega)}^2}\,dx\right)^{\frac{1}{2}}\\
& = \kappa \left(\int_{\Omega} \frac{|u(x,t)|^2}{\|u(\cdot,t)\|_{L^2(\Omega)}^2}\,dx\right)^{\frac{1}{2}}
= {\kappa}.
\end{aligned}
\]
The last inequality contradicts the
assumption that $\,\|d(\cdot,t)\|_{L^2(\Omega)} \le \kappa\,$ in $\,E\,$ unless $\,|\Omega_+(t)\cup \Omega_-(t)| = 0\,$
for almost every $t\in E$. This proves that $\,\|u(\cdot,t)\|_{\Ldue}=0$ almost everywhere in $\,E$. 
With (i) and (ii) proved, the equivalence relation
\eqref{equivalence} is shown. 

The representation \eqref{subdiff} for $\lambda$ follows immediately from the formula for the subdifferential
of $\,g_T$.
\end{proof}

%%%%%%%%%%%%%%%%%%%%%%%%%%%%%%%%%%%%%%%%%%
\subsection{Directional sparsity in time for the optimal control problem}
%%%%%%%%%%%%%%%%%%%%%%%%%%%%%%%%%%%%%%%%%%
The results of the last subsection will now be applied to derive sparsity properties of optimal controls from the variational inequality \eqref{varineq2}. For directional sparsity in time, we use the convex continuous functional
\begin{equation} \label{gT}
g(\bu) =  g(u_1,u_2) : = g_T(u_1) + g_T(u_2) = g_T(I_1\bu) + g_T(I_2\bu),
\end{equation}
where $I_i$ denotes the linear {and continuous projection} mapping $I_i: \bu=(u_1,u_2)  \mapsto u_i$, $i = 1,2$, from $L^1(0,T;L^2(\Omega))^2$ to $L^1(0,T;L^2(\Omega))$.

Since the convex functional $g_T$ is continuous on the whole space $L^1(0,T;L^2(\Omega))$,
we obtain from the sum {and chain rules for subdifferentials (see, e.g., 
\cite[Sect.~4.2.2,~Thm.~1~and\linebreak Thm.2]{ioffe_tikhomirov1979}) that} 
\[
\partial g(\bu) = I_1^*\, \partial g_T(I_1\bu) + I_2^*\, \partial g_T(I_2\bu) = (I,0)^\top \partial g_T(u_1) + (0,I)^\top \partial g_T(u_2),
\] 
with the identical mapping $I \in \mathcal{L}(L^1(0,T;L^2(\Omega)))$. Therefore, we have
\[
\partial g(\bu) = \{ (\lambda_1,\lambda_2) \in L^\infty(0,T;L^2(\Omega))^2: \lambda_i \in \partial g_T(u_i), \, i = 1,2\}.
\]
The variational inequality \eqref{varineq2} is equivalent to two independent variational inequalities for $\overline {u}_1$ and $\overline {u}_2$ that have to hold jointly, namely,
\begin{eqnarray}
\int_Q  \left( - \overline{\mathbf \psi}_1 h({\overline \varphi}) + \kappa \overline \lambda_1 +\nu\, \overline u_1\right)\left(u - \overline{u}_1\right)dx\,dt &\!\!\ge\!\!& 0 \quad \forall\, u \in C_1, \label{varin1}\\
\int_Q  \left( \overline\psi_3
+ \kappa \overline \lambda_2 +\nu\,\overline u_2\right)\left(u - \overline{u}_2\right)dx\,dt 
&\!\!\ge\!\!& 0 \quad \forall\, u \in C_2, \label{varin2}
\end{eqnarray}
where the sets $\,C_i$, $i = 1,2$, are defined by
\[
C_i = \{u \in L^\infty(Q): \underline u_i(x,t)  \le u(x,t) \le \hat u_i(x,t) \mbox{ for a.a. } (x,t) \in Q\},
\]
and where $\overline \lambda_i $, $i = 1,2$, obey {for almost every $t\in (0,T)$} the conditions
\begin{equation} \label{subdiff2}
\overline \lambda_i(\cdot,t) \left\{
\begin{array}{lcl}
\in B(0,1)&\,\mbox{ if } \, \|\overline u_i(\cdot,t) \|_{L^2(\Omega)} = 0\\[1ex]
= \displaystyle \frac{\overline u_i(\cdot,t)}{\|\overline u_i(\cdot,t) \|_{L^2(\Omega)}}&\,\mbox{ if } \, \|\overline u_i(\cdot,t) \|_{L^2(\Omega)} \not= 0
\end{array}
\right. .
\end{equation}
Applying Lemma \ref{L4.1} to each of the variational inequalities \eqref{varin1} and \eqref{varin2} separately, we arrive at the following result:
\begin{theorem}  {\rm (Directional sparsity in time)} \,\,\, Suppose that the general assumptions {\bf (F1)}--{\bf (F4)}
and {\bf (A1)}--{\bf (A6)} are fulfilled, and assume in addition that $\underline u_i,\widehat u_i $ are constants satisfying
\,$\underline u_i<0<\widehat u_i$, for \,$i=1,2$.
Let $\overline \bu = (\overline u_1,\overline u_2)$ be an optimal control of the problem {{\rm (}${\cal CP}${\rm )}} with sparsity functional $\,g\,$ defined in \eqref{gT}, and with associated state $(\overline\mu,
\overline\phi,\overline\sigma)=\S(\ubar)$ solving \eqref{ss1}--\eqref{ss5} and adjoint state $\overline {\boldsymbol{\psi}} = (\overline \psi_1,\overline \psi_2,\overline \psi_3)$ solving  \eqref{adj1}--\eqref{adj5}. Then, there are functions $\overline \lambda_i$, $i=1,2,$ that satisfy \eqref{subdiff2} and \eqref{varin1}--\eqref{varin2}.
In addition, for
almost every $t\in (0,T)$, we have that
\begin{eqnarray}
\|\overline u_1(\cdot,t)\|_{L^2(\Omega)} = 0 \quad &\Longleftrightarrow& \quad \|\overline\psi_1(\cdot,t) h({\overline \varphi}(\cdot,t))\|_{L^2(\Omega)} \le \kappa, \label{u1sparsity}\\
%\end{equation}
%\begin{equation}
\|\overline u_2(\cdot,t)\|_{L^2(\Omega)} = 0 \quad &\Longleftrightarrow& \quad \| \overline\psi_3(\cdot,t)\|_{L^2(\Omega)} \le \kappa. \label{u2sparsity}
\end{eqnarray}
{Moreover, if\, $\overline {\boldsymbol{\psi}}$ and $\overline \lambda_1, \overline \lambda_2$ are given, then}
the optimal controls
$\overline u_1$,  $\overline u_2$ are obtained from the projection formulas
\begin{eqnarray*}
\overline u_1(x,t)& =& \mathbb{P}_{[\underline u_1(x,t) , \widehat u_1(x,t)]} \left(-{\nu}^{-1} \left(-\overline\psi_1(x,t) h({\overline \varphi}(x,t))+ \kappa \overline \lambda_1(x,t)\right)\right),\\
\overline u_2(x,t) &=& \mathbb{P}_{[\underline u_2(x,t) , \widehat u_2(x,t)]} \left(-{\nu}^{-1} \left(\overline\psi_3(x,t) + \kappa \overline \lambda_2(x,t)\right)\right),\quad\mbox{for a.e. $\,(x,t)\in Q$.}
\end{eqnarray*}
\end{theorem}
\Brem
In the medical context, where the controls $\,u_1,u_2\,$ have the meaning of medications or of
nutrients supplied to the patients, it does not seem
to be meaningful to allow for negative controls, unfortunately.
\Erem

{It is to be expected that the support of optimal controls will shrink with increasing sparsity   parameter $\kappa$.
Although this can hardly be quantified or proved, it is useful to confirm that optimal controls vanish for all sufficiently large values of $\kappa$. We are going to derive a corresponding result now.
 
For this purpose, let us indicate for a while the dependence of optimal controls, optimal states, and the associated
adjoint states, on $\kappa$ by an index $\kappa$. An inspection of the conditions \eqref{u1sparsity} and/or \eqref{u2sparsity} reveals that 
$\overline u_{1,\kappa} = 0 $ holds true for all
$\kappa > \kappa_1$, if 
\begin{equation} \label{generalbound1}
\kappa_1 := \sup_{\kappa > 0} \sup_{t\in (0,T)} \| \overline\psi_{1,\kappa}(\cdot,t)h(\overline \varphi_\kappa(\cdot,t))\|_{L^2(\Omega)} < \infty,
\end{equation}
and $\overline u_{2,\kappa} = 0 $ holds true for all
$\kappa > \kappa_2$, if 
\begin{equation} \label{generalbound2}
\kappa_2 = \sup_{\kappa > 0} \sup_{t\in (0,T)} \| \overline\psi_{3,\kappa}(\cdot,t)\|_{L^2(\Omega)} < \infty.
\end{equation}
These boundedness conditions hold simultaneously for $\kappa > \kappa_0=\max\{\kappa_1,\kappa_2\}$. The existence of such a constant $\kappa_0$ will be confirmed next. In order to avoid an overloaded notation, we omit the index $\kappa$ in the following. 

First, we derive bounds for the adjoint state variables $\,\overline\psi_1, \overline\psi_3$
(the function $h(\overline\phi)$ is globally bounded by {\bf (A2)}). To this end, recall the global
estimates \eqref{ssbound1}--\eqref{ssbound3} from Theorem 2.1, which have to be satisfied by 
all possible states $(\mu,\varphi,\sigma)$ corresponding to controls $\bu\in\uad$. It follows that also  
the ``right-hand sides''
$\beta_1 (\overline \varphi - \varphi_Q)$ and $\beta_2 (\overline \varphi(T) - \varphi_\Omega)$ are uniformly bounded, independently of $\kappa$. It remains to show that this implies the boundedness of all possible
adjoint states.

 To this end, recall that by virtue of \eqref{reg1},
\eqref{reg23} we know that $\overline\psi_1\in C^0([0,T];V)$ and $\overline \psi_3\in C^0([0,T];H)$. Now indeed, a
closer look at the proof of \cite[Thm.~2.8]{CSS1} reveals that the bounds derived there are in fact uniform with respect
to the choice of $\bu\in\uad$. Therefore, there is some $\kappa_0>0$ such that $\ubar_{{\kappa}}=\mathbf{0}$ for every 
$\kappa\ge\kappa_0$. For the reader's convenience, we now give some insight how such bounds can be derived.}

In the following, we argue formally, noting that in a rigorous
proof the following arguments would have to be carried out on a Faedo--Galerkin system
approximating the weak form of the adjoint system 
\eqref{adj1}--\eqref{adj5} satisfied by the adjoint variables $(\psi_1,\psi_2,\psi_3)=(\opi_1,\opi_2,\opi_3)$. 
The arguments are similar to those in the proof  of Lemma \ref{L2.4}.

Indeed, we (formally) multiply \eqref{adj1} by $\,\,-\beta\dt\opi_1$, \eqref{adj2} by $\,\,\opi_2$, and
\eqref{adj3} by $\,\,\delta\opi_3$, where $\delta>0$ is yet to be specified. Then we add the three resulting equations, 
whence a cancellation of two terms occurs, and integrate the result over $Q^t:=\Omega\times (t,T)$, where $t\in [0,T)$. 
Using formal integration by parts and the endpoint conditions, we   
then obtain the identity
\begin{align}\label{Elvis}
&{\alpha\beta\int_{Q^t}|\dt\opi_1|^2\,+\,\frac\beta 2\|\nabla\opi_1(t)\|_2^2\,+\frac\beta 2\|\opi_2(t)\|_2^2
\,+\,\frac\delta 2\|\opi_3(t)\|_2^2}\nonumber\\
&{\quad +\int_{Q^t}\bigl(|\nabla\opi_2|^2\,+\,\delta|\nabla\opi_3|^2\bigr)\,+\int_{Q^t}F''_1(\overline\phi)|\opi_2|^2
}\nonumber\\
&{=\,\frac{\beta_2^2}{2\beta}\iO|\overline\phi(T)-\phi_\Omega|^2\,+\,\chi\int_{Q^t}\nabla\opi_2\cdot\nabla\opi_3
\,+\,\beta_1\int_{Q^t}(\overline\phi-\phi_Q)\opi_2\,-\int_{Q^t}F_2''(\overline\phi)|\opi_2|^2}\nonumber\\
&{\qquad +\int_{Q^t}(P\overline\sigma-A-\overline u_1)h'(\overline\phi)\opi_1 \opi_3\,-\int_{Q^t}
E\overline\sigma h'(\overline\phi)\opi_2\opi_3\,-\,\delta\int_{Q^t}(Eh(\overline\phi)+B)|\opi_3|^2}\nonumber\\
&{\qquad+\,\delta\int_{Q^t}(P h(\overline\phi)\opi_1+\chi\opi_2)\opi_3. }
\end{align}
{Since $F_1''\ge 0$, all of the terms on the left-hand side are nonnegative. Moreover, Young's inequality implies
that}
$${\chi\int_{Q^t}\nabla\opi_2\cdot\nabla\opi_3\,\le\,\frac 12\int_{Q^t}|\nabla\opi_2|^2\,+\,\frac{\chi^2}2
\int_{Q^t}|\nabla\opi_3|^2.}$$
{Hence, invoking the known bounds for the state variables, and applying Young's inequality appropriately to 
the terms on the right-hand side, we obtain from \eqref{Elvis} the estimate }
\begin{align}
\label{Mick}
&{\alpha\beta\int_{Q^t}|\dt\opi_1|^2\,+\,\frac\beta 2\|\nabla\opi_1(t)\|_2^2\,+\frac\beta 2\|\opi_2(t)\|_2^2
\,+\,\frac\delta 2\|\opi_3(t)\|_2^2}\nonumber\\
&{\quad +\,{\frac 12}\int_{Q^t}|\nabla\opi_2|^2\,+\,\left(\delta-\mbox{$\frac 12$}\,\chi^2\right)\int_{Q^t}|\nabla\opi_3|^2
}\nonumber\\
&{\le \,C_1\,+\,C_2(1+\delta)\int_{Q^t}\left(|\opi_1|^2+|\opi_2|^2+|\opi_3|^2\right),}
\end{align}
with constants $C_1,C_2$ that depend neither on $\uad$ nor on $\,\kappa$. 

Next observe that $\,\opi_1(T)=0\,$ and thus $\,\,\frac 12\,\|\opi_1(t)\|_2^2=-\int_t^T(\dt\opi_1(s),\opi_1(s))\,ds$. Hence,
owing to Young's inequality,
\begin{equation}
\label{Jagger}
\frac 12\,\|\opi_1(t)\|_2^2\,\le\,\frac{\alpha\beta}2\int_{Q^t}|\dt\opi_1|^2\,+\,\frac 1{2\alpha\beta}
\int_{Q^t}|\opi_1|^2\,.
\end{equation}
Now we add \eqref{Mick} and \eqref{Jagger} and choose $\,\delta=\chi^2$. Using Gronwall's lemma 
backward in time, it then easily follows that, in particular,
$$
{\|\opi_1\|_{L^\infty(0,T;V)}\,+\,\|\opi_3\|_{L^\infty(0,T;H)}\,\le\,C_3,}
$$
 where $C_3>0$ is independent of both $\uad$ and $\,\kappa$. Then,
$$
{\|\opi_1 h(\overline\phi)\|_{L^\infty(0,T;H)}\,+\,\|\opi_3\|_{L^\infty(0,T;H)}\,\le\,\left(1+\|h\|_{L^\infty
(\erre)}\right)C_3}\,=:\,\kappa_0.
$$
{The asserted existence of the constant $\,\kappa_0\,$ is thus shown.}

%%%%%%%%%%%%%%%%%%%%%%%%%%%%%%%%%%%%%%%%%%
\subsection{Spatial directional sparsity and full sparsity}
%%%%%%%%%%%%%%%%%%%%%%%%%%%%%%%%%%%%%%%%%%

Let us briefly sketch the other types of sparsity that are obtained {from the choices} $\,g = g_\Omega\,$ and
$\,g = g_Q$, respectively.

With the functional $\,g_\Omega$, we obtain regions in $\Omega$ where the optimal controls are zero 
for a.e. $t\in (0,T)$. The theory is analogous to that of directional sparsity in time: indeed, it is obtained by
simply interchanging the roles of $t$ and $x$. For instance, instead of the equivalences \eqref{u1sparsity}, \eqref{u2sparsity},
one obtains for a.e.  $x \in \Omega$ that
\begin{eqnarray*}
\|\overline u_1(x,\cdot)\|_{L^2(0,T)} = 0 \quad &\Longleftrightarrow& \quad \|\overline\psi_1(x,\cdot) h({
\overline \varphi}(x,\cdot))\|_{L^2(0,T)} \le \kappa, \\
%\end{equation}
%\begin{equation}
\|\overline u_2(x,\cdot)\|_{L^2(0,T)} = 0 \quad &\Longleftrightarrow& \quad \| \overline\psi_3(x,\cdot)\|_{L^2(0,T)} \le \kappa. 
\end{eqnarray*}
For the choice $\,g = g_Q$, the equivalence relations 
\begin{eqnarray*}
\overline u_1(x,t) = 0 \quad &\Longleftrightarrow& \quad |\overline\psi_1(x,t) h({\overline \varphi}(x,t))| 
\le \kappa, \\
%\end{equation}
%\begin{equation}
\overline u_2(x,t) = 0 \quad &\Longleftrightarrow& \quad | \overline\psi_3(x,t)| \le \kappa,
\end{eqnarray*}
can be deduced for almost every $\,(x,t) \in Q$. We refer to the discussion of the variational inequality \eqref{varineqaux} in \cite{casas_ryll_troeltzsch2013}. Therefore, the optimal controls vanish in certain spatio-temporal subsets of $\,Q$.

Moreover, in this case a usually unexpected property of the function $\overline {\boldsymbol \lambda } \in g(\overline \bu)$ is obtained: $\overline {\boldsymbol \lambda }$ is unique, that is, for an optimal control, the subdifferential is a
 singleton; we again refer to \cite{casas_ryll_troeltzsch2013}. This fact can easily be explained. Consider, e.g., the function $\overline \lambda_2 \in \partial g_Q(\overline u_2)$:

Thanks to \eqref{dgQ}, it holds that
\[
\overline \lambda_2(x,t) = \left\{
\begin{array}{rcl}
1& \mbox{ if } &\overline u_2(x,t) > 0\\
-1& \mbox{ if } &\overline u_2(x,t) < 0
\end{array}
\right.
\]
Therefore, the only points, at which $\overline \lambda_2(x,t)$ might not be uniquely determined, are those where $\overline
u_2(x,t)$ vanishes. At these points, however, $\overline u_2(x,t)=0\,$ is away from the {thresholds}, and hence the reduced gradient must be zero, i.e.,
\[
0 = \overline\psi_3(x,t) + \kappa \overline \lambda_2(x,t) + \nu \cdot 0.
\]
This implies that $\overline \lambda_2(x,t) = - \kappa^{-1} \overline\psi_3(x,t)$ {at} these points. With a little more effort, finally the projection formula
\[
\overline \lambda_2(x,t) = \mathbb{P}_{[-1,1]}\left(- \frac{1}{\kappa} \overline\psi_3(x,t)\right)
\]
results. {By similar reasoning, the identity}
\[
\overline \lambda_1(x,t) = \mathbb{P}_{[-1,1]}\left( \frac{1}{\kappa}  
\overline\psi_1(x,t) h({\overline \varphi}(x,t)) \right)
\]
can be derived.

%%%%%%%%%%%%%%%%%%%%%%%%%%%%%%%%%%%%%%%%%%%%%%%%%%%%%%%%%%%%%%%%%

%\bibliographystyle{plain}
%\bibliography{fredi.bib}
%
%\end{document}
%
%
%
%\vspace{3truemm}

\Begin{thebibliography}{10}

\bibitem{cartan1967}
H.~Cartan:
``Calcul diff\'erentiel. Formes diff\'erentielles''.
Hermann, Paris, 1967.

\bibitem{casas_herzog_wachsmuth2017}
E.~Casas, R.~Herzog, G.~Wachsmuth:
Analysis of spatio-temporally sparse optimal control problems of
  semilinear parabolic equations.
{\em ESAIM Control Optim. Calc. Var.} {\bf 23} (2017), 263--295.

\bibitem{casas_ryll_troeltzsch2013}
E.~Casas, C.~Ryll, F.~Tr\"{o}ltzsch:
Sparse optimal control of the {S}chl\"{o}gl and {F}itz{H}ugh--{N}agumo
systems.
{\em Comput. Methods Appl. Math.} {\bf 13} (2013), 415--442.

\bibitem{casas_ryll_troeltzsch2015}
E.~Casas, C.~Ryll, F.~Tr\"{o}ltzsch:
Second order and stability analysis for optimal sparse control of the
  {F}itz{H}ugh-{N}agumo equation.
 {\em SIAM J. Control Optim.} {\bf  53} (2015), 2168--2202.

%\bibitem{casas_troeltzsch2019}
%E.~Casas and F.~Tr\"{o}ltzsch.
%\newblock Optimal sparse boundary control for a semilinear parabolic equation
%  with mixed control-state constraints.
%\newblock {\em Control and Cybernetics}, 48(1):89--124, 2019.

\bibitem{CGH}
P. Colli, G. Gilardi, D. Hilhorst:
On a Cahn--Hilliard type phase field system related to tumor growth. {\em Discret. Cont. Dyn. Syst.} {\bf 35}
(2015), 2423--2442.                                          

\bibitem{CGRS1}
P. Colli, G. Gilardi, E. Rocca, J. Sprekels: Vanishing viscosities and error estimate for a Cahn--Hilliard
type phase field system related to tumor growth. {\em Nonlinear Anal. Real World Appl.} {\bf 26} (2015),
93--108.

\bibitem{CGRS2}
P. Colli, G. Gilardi, E. Rocca, J. Sprekels: Asymptotic analyses and   error estimates for a Cahn--Hilliard
type phase field system modelling tumor growth. {\em Discret. Contin. Dyn. Syst. Ser. S} {\bf 10} (2017),
37--54.

\bibitem{CGRS3}
P. Colli, G. Gilardi, E. Rocca, J. Sprekels:
 Optimal distributed control of a diffuse interface model of tumor growth. {\em
Nonlinearity} {\bf 30} (2017), 2518--2546.

\bibitem{CGS24}
P. Colli, G. Gilardi, J. Sprekels:
A distributed control problem for a fractional tumor growth model. {\em Mathematics} {\bf 7} (2019), 792.

\bibitem{CSS1}
P. Colli, A. Signori, J. Sprekels: Optimal control of a phase field system
modelling tumor growth with chemotaxis and singular potentials. {\em Appl. Math. Optim.},
Online First, October 21, 2019, https://doi.org/10.1007/s00245-019-09618-6.

%\bibitem{DHP}
%R. Denk, M. Hieber, J. Pr\"uss: Optimal $L^p-L^q-$ estimates for parabolic boundary value 
%problems with inhomogeneous data. {\em Math. Z.} {\bf 257} (2007), 193--224.

\bibitem{Dieu}
J. Dieudonn\'e: ``Foundations of Modern Analysis''. Pure and Applied Mathematics, vol. 10,
Academic Press, New York, 1960.

\bibitem{GLSS}
H. Garcke, K.\,F. Lam, E. Sitka, V. Styles: A Cahn--Hilliard--Darcy model for tumour
growth with chemotaxis and active transport. {\em Math. Model. Methods Appl. Sci.} {\bf 26} (2016),
1095--1148.

\neuj{
\bibitem{GT}
D. Gilbarg, N.\,S. Trudinger: ``Elliptic partial differential equations of second order'' (2nd edition),
Springer-Verlag, Berlin-Heidelberg, 1983.
}

\bibitem{HZO}
A. Hawkins-Daarud, K.\,G. van der Zee, J.\,T. Oden: Numerical simulation of a thermodynamically
consistent four-species tumor growth model. {\em Int. J. Numer. Math. Biomed. Eng.} {\bf 28} (2011),
3--24.

\bibitem{ekeland_temam1974}
I.~Ekeland and R.~Temam:
``Analyse convexe et probl\`emes variationnels''.
Dunod, Gauthier-Villars, Paris-Brussels-Montreal, Que., 1974.

\bibitem{herzog_obermeier_wachsmuth2015}
R.~Herzog, J.~Obermeier, G.~Wachsmuth:
Annular and sectorial sparsity in optimal control of elliptic
  equations.
{\em Comput. Optim. Appl.} {\bf 62} (2015), 157--180.

\bibitem{herzog_stadler_wachsmuth2012}
R.~Herzog, G.~Stadler, G.~Wachsmuth:
Directional sparsity in optimal control of partial differential
  equations.
{\em SIAM J. Control Optim.} {\bf 50} (2012), 943--963.       

\bibitem{ioffe_tikhomirov1979}
A. D. Ioffe, V. M. Tikhomirov: ``Theory of extremal problems''. Studies in
  Mathematics and its Applications, vol. 6, North-Holland Publishing Co., Amsterdam-New York, 1979.

%\bibitem{LSU}
%O.\,A. Lady\v{z}enskaja, V.\,A. Solonnikov, N.\,N. Uralceva: ``Linear and Quasilinear Equations of Parabolic
%Type''. Mathematical Monographs, vol. 23, American Mathematical Society, Providence, Rhode Island, 1968.

\bibitem{LM}
J.\,L. Lions, E. Magenes: ``Non-Homogeneous Boundary Value Problems'', vol. I, Springer-Verlag, Heidelberg, 1972.

\bibitem{Simon}
J. Simon: Compact sets in the space $L^p(0,T;B)$. {\em Ann. Mat. Pura Appl.} {\bf 146} (1997), 65--96.

\bibitem{stadler2009}
G.~Stadler:
Elliptic optimal control problems with {$L^1$}-control cost and
  applications for the placement of control devices.
{\em Comput. Optim. Appl.} {\bf 44} (2009), 159--181.

\bibitem{Fredibuch}
F. Tr\"oltzsch: ``Optimal Control of Partial Differential Equations: Theory, Methods and Applications''.
Graduate Studies in Mathematics vol. 112, American Mathematical Society, Providence, Rhode Island, 2010.

\End{thebibliography}

\End{document}

%%%%%%%%%%%%%%%%%%%%%%%%%%%%%%%%%%%%%%%%%%%%%